\documentclass[10pt,a4paper]{article}
\usepackage[utf8]{inputenc}
\usepackage{amsmath}
\usepackage{yfonts}
\usepackage{amsthm}
\usepackage{amsfonts}
\usepackage{amssymb}
\usepackage{geometry}       
\usepackage{pgfplots}
\pgfplotsset{compat=1.18}
\usepackage{caption}
\usepackage{graphicx}	
\usepackage{tikz}
\usepackage{subfigure}
\usepackage{float}	
\usepackage{colortbl}
\usepackage{xcolor}
\usepackage{comment}
\usepackage{hyperref}
\usepackage{soul,color}
\usepackage{tabularx}
\usepackage{longtable}

\theoremstyle{definition}
\newtheorem*{thm*}{Theorem}
\newtheorem{defi}{Definition}[section]

\theoremstyle{remark}

\theoremstyle{plain}
\newtheorem*{fact}{Fact}
\newtheorem{lem}[defi]{Lemma}
\newtheorem{thm}[defi]{Theorem}

\newtheorem{rem}[defi]{Remark}
\newtheorem{example}{Example}
\author{F. Ciavattini \footnote{Doctoral School in Mathematics, University of Rome ``Tor Vergata'' (Italy)}, 
A. Della Corte\footnote{Mathematics Division, School of Sciences and Technology, University of Camerino (Italy). 
\\ Corresponding author - email to: alessandro.dellacorte@unicam.it.}, 
C. Lucamarini \footnote{Doctoral school in Ergodic Theory of Dynamical Systems, Università degli Studi di Pisa (Italy).}
}
\title{Order spectra in topological dynamical systems}
\begin{document}

\maketitle

\begin{abstract}
\noindent In a compact topological dynamical system $(X,f)$, we associate to every pair $(x,y)$ a canonical order-theoretic invariant, its emergent order spectrum $\Omega(x,y)$. We first prove that, if $x$ and $y$ are chain-related, one can always build families of nested and acyclic $\varepsilon_n$-chains ($\varepsilon_n \to 0$). The order spectrum $\Omega(x,y)$ is then defined as the set of countable linear order-types obtained as direct limits of (order-compatible) nested and acyclic $\varepsilon_n$-chains. The spectrum is independent of the vanishing sequence and invariant under topological conjugacy, and it refines Conley’s partial order.


\end{abstract}

\section{Introduction}

Many questions in dynamics, from structural stability to the behavior of metastable systems, concern what can reliably be said about the structure of dynamics when one allows small perturbations. One way to formalize this is to replace exact orbits by $\varepsilon$–chains (or pseudo-orbits) and study the relations they generate. The notion of an $\varepsilon$–chain thus plays a central role in dynamics, providing the foundation for both  C.~Conley’s decomposition theory of dynamical systems \cite{co78} and E.~Akin’s structural description of attractors (see, e.g., \cite{kurka03}, Theorem~2.68, p.~82).

A point $x$ is chain-related to a point $y$ (and we write $x\,\mathcal{C}\,y$) if, for every $\varepsilon>0$, it is possible to go from $x$ to $y$ by repeating finitely many times the operation of ``taking the $f$-image and applying a correction smaller than $\varepsilon$''. A finite sequence of points
$$x_0=x,\,x_1,\dots,\,x_n=y$$
such that $d(f(x_i),x_{i+1})<\varepsilon$, for every $i\in\{0,\dots,n-1\}$, is called an $\varepsilon$-chain.

The starting point of this work is the observation that the chain relation hides much more structure than the mere existence of arbitrarily fine pseudo-orbits. Between two chain-related points $x,y$ there are in general many $\varepsilon$–chains, and each of these chains comes equipped with a natural linear order (given by the indices along the chain). At first sight, the limit 
$\varepsilon\to 0$ seems to remember only whether two points belong to the chain relation. The definition itself gives no immediate indication that the finite-$\varepsilon$ chains should leave behind any structurally invariant “track”, and one might even suspect that no such track exists at all. Herein we show that it is indeed possible to identify a canonical track encoded in the invariant structure of the chains. 

This latent structure can be made explicit once realized that it is always possible to pass to families of chains that are i) nested and ii) acyclic. These chains provide a structure supporting linear orders that become increasingly representative of the recurrence properties as $\varepsilon_n$ decreases. Then we record only the order-types that are limits of (order-compatible) sequences of nested, acyclic chains. 
This produces, for each ordered pair $(x,y)$ of points in the space a (possibly empty) collection of countable linear order-types, the \emph{emergent order spectrum} (EOS), denoted by $\Omega(x,y)$.

Our main result is the following:
\begin{thm*}
\emph{Let $(X,f)$ be a topological dynamical system, with $X$ a compact metric space and $f:X\longrightarrow X$ continuous. Given any pair of chain-related points $x,y$, there is a countable family $\{C_n\}$ of $\varepsilon_n$–chains (with $\varepsilon_n \to 0$) such that:
\begin{enumerate} 
\item for every $n$, $C_n\subseteq C_{n+1}$; 
\item every $\epsilon_n$-chain is acyclic, that is it does not contain any cyclic sub-chain.
\end{enumerate}
Any order-compatible sequence of nested, acyclic $\varepsilon_n$-chains produces a unique limit order-type supported on $\bigcup_n C_n$. The collection of the order-types so obtained for the pair $(x,y)$ is independent of the sequences $\{\varepsilon_n\}$ and of the compatible metric, and invariant under conjugacy.}
\end{thm*}

This result is proven in three steps. We prove point 1. in Theorem~\ref{C=C_nested}, through a projection on the Hausdorff limit of a Hausdorff-convergent sequence of $\varepsilon_n$-chains. The existence of nested chains for chain-related points is false, in general, in non-compact spaces (an example is provided). Point 2. is proven in Theorem~\ref{mainth} through a transfinite iteration of the Hausdorff projection arriving at a stabilizing closed, invariant limit supporting nested and acyclic chains. 

Finally, to prove the last statement, we extract order-theoretic information from  order-compatible nested, acyclic chains $\{C_n\}$. This leads to the notion of EOS introduced above, which captures the types of recurrence exhibited. 

\noindent In the last section, we show that the EOS contains quite comprehensive
information on the invariant recurrence structure of the system.
A structural property of the EOS is that it refines Conley’s
decomposition: for each pair of chain-related points $(x,y)$ and each
$\beta\in\Omega(x,y)$, the corresponding limit linear order on the union
of the chains decomposes into convex blocks indexed by the chain components
encountered by the chains, and the induced order on these blocks is precisely
the Conley partial order restricted to those components
(Theorems~\ref{convex_} and~\ref{conley_}). Moreover, we provide examples
in which the EOS refines Auslander's prolongational sets, distinguishing
recurrence patterns in cases in which prolongational ranks coincide.


\section{Preliminaries}
By $\mathbb{N}$ and $\mathbb N_0$ we mean, respectively, the set of positive integers and of non-negative integers,
while $\mathbb Z$, $\mathbb{Q}$ and $\mathbb R$ indicate, respectively, the set of integers, rational and real numbers. We denote by $\mathbb S^1$ the standard unit circle.

By compact system we mean a pair $(X,f)$, where $X$ is a metric space, with metric $d$, and
$$f:X\longrightarrow X$$
is a continuous map. We say that the system $(X,f)$ is compact if $X$ is a compact space. Given a point $x\in X$, we denote by $\mathcal{O}(x)$ its orbit, that is the set
$$\mathcal{O}(x)=\{f(x),f^2(x),\dots\}.$$

For $S\subseteq X$, we write $\partial S$ for the topological boundary of $S$, and $\operatorname{int}(S)$ for the topological interior of $S$, respectively, and we denote by $\overline{S}$ the topological closure of $S$. We indicate the open (closed) ball of radius $\rho$ centered at $x$ by $B_\rho(x)$ ($\overline{B}_\rho(x)$).

Let us review the well-known topological dynamical relations $\mathcal{O}\subseteq\mathcal{R}\subseteq\mathcal{N}\subseteq\mathcal{C}\subseteq X^2$, as introduced by E. Akin in \cite{Akin} for the more general case of closed relations of $X$. In the case of a continuous map, they can be defined as follows (see for instance [\cite{kurka03}, Def. 2.2, p. 47]).

\begin{defi}   
\mbox{}
We introduce the following standard topological dynamical relations:

    \begin{itemize}
        \item[1)] Orbit relation:\\ $x\,\mathcal{O}\,y$ if and only if $\exists k\in\mathbb{N}$ such that $f^k(x)=y$;
        \item[2)] Recurrence relation:\\$x\,\mathcal{R}\,y$ if and only if $\forall\,\varepsilon>0\, \exists k\in\mathbb{N}$ such that $f^k(x)\in B_\varepsilon (y)$;
        \item[3)] Non-wandering relation:\\ $x\,\mathcal{N}\,y$ if and only if $\forall\,\varepsilon>0\, \exists z\in B_\varepsilon (x)$ and $\exists k\in\mathbb{N}$ such that $f^k(z)\in B_\varepsilon (y)$;
    \end{itemize}
\end{defi}

Let us also recall the notion of $\varepsilon$-chain and chain relation.
\begin{defi}\label{chain}
    Let $(X, f)$ be a dynamical system. Given two points $x, y \in X$ and $\varepsilon > 0$, an $\varepsilon$-chain from $x$ to $y$ is an indexed, finite sequence of points of $X$, that is a map
    
    $$C:\{0,1,\dots,m\}\longrightarrow X,$$
    
    with $m\ge1$, such that, setting $x_i=C(i)$, we have:
    \begin{itemize}
        \item[i)] $x_0 = x$ and $x_m = y$,
        \item[ii)] $d(f(x_i), x_{i+1})<\varepsilon$, for every $i=0,1,\dots,m-1$.
    \end{itemize}
    By $\widehat{C}$, we denote the support of the chain $C$, that is, 
    $$\widehat{C}=C(\{0,1,\dots,m\})=\{x_0,\dots,x_m\}.$$ 
    
    With a (very common) abuse of notation, we may indicate a chain $C$ by writing simply $$C:x_0,\,x_1,\dots,\,x_m.$$

    We say $x,y \in X$ are in \emph{chain recurrence relation} (or simply in chain relation), and we write $x\,\mathcal{C}\,y$, if and only if, for every $\varepsilon>0$, there exists an $\varepsilon$-chain from $x$ to $y$.
\end{defi}

The relation $x\equiv_C y$ defined as $x\equiv_C y$ if and only if $x\,\mathcal{C}\, y$ and $y\,\mathcal{C}\, x$ is an equivalence relation. Its equivalence classes are called \textit{chain components}. The chain component containing $x$ will be denoted by $[x]$.

\begin{defi}\label{compchain}
    Let $(X,f)$ be a dynamical system.
    Assuming that, for $x,y\in X$, we have $x\,\mathcal{C}\,y$, we say that the sequence $\{C_n\}_{n=1}^\infty$ of $\varepsilon_n$-chains from $x$ to $y$ is a \emph{complete sequence of chains} if $\varepsilon_n$ converges to 0 monotonically.
\end{defi}

\begin{defi}
    Let $(X,f)$ be a compact dynamical system and $K,K'\subseteq X$ be two chain components. We set $$K\le_{\mathrm{Conley}}K'\iff \ \exists x \in K\,,\,y\in K'\,:\, x\,\mathcal{C}\,y.$$ 
    The relation $\le_{\mathrm{Conley}}$ is a partial order on the chain components.
\end{defi}

\begin{defi}
    Given a linear order $(S,\preceq)$, we say that $M\subseteq S$ is \textit{convex} with respect to $\preceq$ (or simply convex if no confusion may arise) if, for every $x,y\in M$, we have $x\preceq z\preceq y\implies z\in M$. 
\end{defi}

\begin{defi}
    Let $(X, f)$ be a dynamical system. We say that $(X,f)$ is \emph{transitive} if and only if there exists $x\in X$ such that $\overline{\mathcal{O}(x)}=X$.
\end{defi}

Let us also recall some standard definitions of basic concepts in topological dynamics (see for instance \cite{kurka03}). 

\begin{defi}
    Let $(X,f)$ be a compact dynamical system. 
\begin{itemize}
    \item  The $\omega$-limit set of $x\in X$ is defined as:
    $$\omega(x) := \bigcap_{N\geq 0} \overline{\{f^k(x) : k\ge N\}}.$$ 

    \item A closed set $U\subseteq X$ is
    called an \emph{inward set} if $f(U)\subseteq\operatorname{int}(U)$. 

    \item  An \emph{attractor} is a set $A\subseteq X$ for which there exists 
    an inward set $U$ such that
    $$A = \bigcap_{n\geq 0} f^n(U).$$

    \item The \emph{basin} of an attractor $A$ is
    $$B(A) := \{x\in X : \omega(x)\subseteq A\}.$$

    \item     The \emph{dual repeller} of an attractor $A$ is
    $$R := X \setminus B(A).$$
\end{itemize}
\end{defi}

\begin{defi}
    Given two closed sets $F_1,F_2\subseteq X$, the Hausdorff distance $d_{\mathcal{H}}(F_1,F_2)$ between them is given by (see for instance \cite{munkres2013}, Section 45):
    $$d_\mathcal{H}(F_1,F_2):=\inf\{\varepsilon>0\,:\,F_1\subseteq B_\varepsilon(F_2)\text{ and }F_2\subseteq B_\varepsilon(F_1)\}.$$
\end{defi}
\begin{rem}
    If $d_\mathcal{H}(F_1,F_2)=\varepsilon$, then $F_1\subseteq\overline{B}_\varepsilon(F_2)$ and $F_2\subseteq\overline{B}_\varepsilon(F_1)$.
    Let $\mathfrak{F}$ be the collection of nonempty closed sets of $X$. It is known that the metric space $(\mathfrak{F},d_\mathcal{H})$ is compact if and only if $X$ is compact. Thus, assuming that $X$ is compact, for every sequence of closed sets $\{F_n\}_{n=1}^\infty$, there is a subsequence $\{F_{n_k}\}_{k=1}^\infty$ and a closed set $F_\infty\subseteq X$ such that $d_\mathcal{H}(F_{n_k},F_\infty)\xrightarrow{k\to\infty}0$.
\end{rem}
\begin{rem}
    We recall that, in a compact space, we always have:
\begin{equation}\label{nest_haus}
A_1\subseteq A_2\subseteq\dots \subseteq A_n\subseteq\dots \implies A_n\xrightarrow{n\to\infty}\overline{\bigcup_{n=1}^\infty A_n}
\end{equation}
where convergence is meant with respect to the Hausdorff distance (see Corollary 5.32 in \cite{Tuzhilin}).
\end{rem}

\begin{defi}
    Let $(P,\leq_P)$ and $(Q,\leq_Q)$ be posets. An \emph{order-isomorphism} between $P$ and $Q$ is a bijection $f : P \to Q$ such that for all $x,y \in P$,
    $$x \leq_P y \quad \Longleftrightarrow \quad f(x) \leq_Q f(y).$$
    We write $(P,\leq_P) \cong (Q,\leq_Q)$ if such a map exists. This defines an equivalence relation on the class of all posets.

    The \emph{order-type} of a poset $(P,\leq_P)$ is the equivalence class
    $$\operatorname{otp}(P,\leq_P) \;=\; \{ (Q,\leq_Q) : (Q,\leq_Q) \cong (P,\leq_P) \}.$$
\end{defi}

\begin{rem}
    Of course the collection $\operatorname{otp}(P,\leq_P)$ is in general not a set 
    but a proper class in ZF. 
    One common way to avoid this size issue is to work inside a 
    \emph{universe of sets} (for example, the cumulative hierarchy $V$ 
    in von Neumann's construction (see for instance \cite{Jech}, p.63)), and to take order-types with respect to 
    the posets contained in that universe. 
\end{rem}

Let us now define two binary relations $\mathcal{C}_\subseteq,\mathcal{C}_\preceq\subseteq X^2$ that are for sure not weaker than the chain relation $\mathcal{C}$.

\begin{defi}\label{nested}
    Let $(X,f)$ be a dynamical system. For $x,y\in X$, we say that $x\,\mathcal{C}_\subseteq\, y$ if and only if, for a sequence $\{\varepsilon_n\}_{n=1}^\infty$ of positive real numbers converging to 0 monotonically, there is a complete sequence $\{C_n\}_{n=1}^\infty$ of $\varepsilon_n$-chains from $x$ to $y$ such that $\widehat{C_n}\subseteq \widehat{C_{n+1}}$.

    We will call $\{C_n\}_{n=1}^\infty$ a \textit{complete sequence of nested chains for the sequence $\{\varepsilon_n\}_{n=1}^\infty$}.
\end{defi}

\begin{rem}\label{canon}
    Def.\ref{nested} is independent of the particular sequence chosen, since, given any sequence $\{\delta_m\}_{m=1}^\infty$ of positive real numbers converging to 0, there exists a sub-sequence $\{C_{n_m}\}_{m=1}^\infty$ of nested $\delta_m$-chains. For this reason, in the following we may consider  collections of nested chains for a certain pair of points in $X$ without specifying the sequence of real numbers.
\end{rem}

\begin{defi}
    Let $(X,f)$ be a dynamical system.
    For $x,y\in X$ and $\varepsilon>0$ let $C:x_0,\,x_1,\dots,\,x_m$ be an $\varepsilon$-chain  from $x=x_0$ to $y=x_m$. We say that $C$ contains a \emph{cyclic sub-chain} if one of the following conditions holds:
    \begin{itemize}
        \item for some $i,j\in\{1,\dots,m-1\}$ such that $i\neq j$, we have $x_i=x_j$;
        \item for some $i\in\{1,\dots,m-1\}$, we have $x_i=x_0$ or $x_i=x_m$.
    \end{itemize}
    We say that $C$ is \emph{acyclic} if it does not contain any cyclic sub-chain.
\end{defi}

    It is easily seen that any $\varepsilon$-chain $C$ with a cyclic sub-chain can be modified (by suitably eliminating some points) to get an $\varepsilon$-chain $C'$ with no cyclic sub-chain and such that $\widehat{C'}\subseteq \widehat{C}$.

\begin{defi}
    Given a sequence of totally ordered sets $\{(A_n,\preceq_n)\}_{n=1}^\infty$ such that $A_n\subseteq A_{n+1}$ for every $n\in\mathbb N$, we say that
    \begin{enumerate}
        \item $\{(A_n,\preceq_n)\}_{n=1}^\infty$ is \emph{weakly order-compatible} if, for every element $a,b\in\cup_nA_n$, there exists $N\in\mathbb N$ such that exactly one of the following holds:
        \begin{itemize}
            \item $a\preceq_nb$, for every $n\ge N$, or,
            \item $b\preceq_na$ for every $n\ge N$;
        \end{itemize}
        \item $\{(A_n,\preceq_n)\}_{n=1}^\infty$ is \emph{order-compatible}, if, for every $n\in\mathbb N$, for every $a,b\in A_n$ such that $a\preceq_nb$, we have $a\preceq_{n+1}b$, too.
    \end{enumerate}
\end{defi}
\begin{rem}\label{diag_}
    Given a sequence of totally ordered sets $\{(A_n,\preceq_n)\}_{n=1}^\infty$ such that $A_n\subseteq A_{n+1}$ for every $n\in\mathbb N$, we can always find a weakly order-compatible sub-sequence $\{(A_{n_k},\preceq_{n_k})\}_{k=1}^\infty$ by  diagonalization. In general, however, there is no order-compatible subsequence.
\end{rem}

\begin{defi}\label{preceq}
    Let $(X,f)$ be a dynamical system.
    Assume that two points $x,y\in X$ are such that $x\,\mathcal{C}_\subseteq\,y$. We say that $x\,\mathcal{C}_\preceq\,y$ if the sequence of nested chains $\{C_n\}_{n=1}^\infty$ can be chosen so as to make all its chains without cyclic sub-chains.
\end{defi}

\begin{rem}\label{canon_}
    Also here, if $x\,\mathcal{C}_\preceq\,y$, then the collection of acyclic nested $\varepsilon$-chains $\{C_n\}_{n=1}^\infty$ works (passing if needed for a subsequence) for every positive sequence $\varepsilon_n$ converging to zero. Thus, the fact that $x$ is in $\mathcal{C}_\preceq$-relation with $y$, is independent of the specific sequence $\{\varepsilon_n\}_{n=1}^\infty$. For this reason, in the following we may consider sequences of acyclic, nested chains for a certain pair of points in $X$ without specifying the sequence of real numbers.
\end{rem}

\begin{rem}\label{lem:chain-indipendence-metric}
 Let $(X,f)$ be a compact dynamical system, and let $d$ be the metric on $X$.
 Consider $x,y\in X$ such that $x\,\mathcal{C}\,y$, and let 
 $\{C_n\}_{n=1}^\infty$ be a complete sequence of chains from $x$ to $y$.
 Then it is easily seen that, for every other metric $d'$ on $X$ that is equivalent to $d$, 
 the sequence $\{C_n\}_{n=1}^\infty$ is again a complete sequence of chains 
 from $x$ to $y$.
\end{rem}

In the following, we will prove that in fact $\mathcal{C}_\subseteq$ and $\mathcal{C}_\preceq$ coincide with $\mathcal{C}$ in compact systems. 

It may seem, at first sight, that the challenging part in proving the equalities $\mathcal{C}=\mathcal{C}_\subseteq=\mathcal{C}_\preceq$ is to build nested chains, while the avoidance of cycles is a technical nuisance one can take care of easily. As we will see when we will address the problem in Section \ref{sec5}, things are different. 

\section{\texorpdfstring{$\mathcal{C}=\mathcal{C}_\subseteq$}{C=C\_subseteq}}\label{section_Cnest}

In this section, we prove that, in the compact case, we can always build nested chains for two points in chain relation (Theorem \ref{C=C_nested}). We will refine later the construction to achieve nested chains with compatible linear orderings.

Let $\{C_k\}_{k=1}^\infty$ be a complete sequence of chains. We want to obtain a sequence of nested chains $\{S_n\}_{n=1}^\infty$ whose supports are contained in the Hausdorff limit of a converging subsequence of $\{\widehat{C_k}\}_{k=1}^\infty$. We will call this sequence a \textit{Hausdorff projection} of the original sequence of chains (Def.\ref{hauspro}).

    Since $\widehat{C_k}$ is a finite set for every $k\in\mathbb N$, $\{\widehat{C_k}\}_{k=1}^\infty$ is a sequence of closed sets. Up to considering subsequences, we can assume that 
    \begin{equation}
    \label{hausconv}
        \widehat{C_k}\xrightarrow{k\to\infty}C_\infty,
    \end{equation}
    where $C_\infty$ is a closed set and the convergence is in the Hausdorff metric.

    Setting $\eta=4\,\cdot\,\mbox{diam}(X)$ and $\varepsilon_n=\frac{\eta}{2^n}$, we will prove that we can find a sequence of nested chains $\{S_n\}_{n=1}^\infty$, where, for every $n$, $S_n$ is an $\varepsilon_n$-chain from $x$ to $y$ such that $\widehat{S_n}\subseteq\widehat{S_{n+1}}\subseteq C_\infty$.
    We proceed by induction on $n\in\mathbb{N}$. 
    Since $x,y\in \widehat{C_n}$ for every $n\in\mathbb N$, we have $x,y\in C_\infty$ and, therefore, $S_1:x_0=x,\,x_1=y$ is an $\varepsilon_1$-chain.

    Assume now that we have an $\varepsilon_n$-chain $S_n:x_0=x,\dots,x_m=y$ between $x$ and $y$ such that $\widehat{S_n}\subseteq C_\infty$. We want to build an $\varepsilon_{n+1}$-chain $S_{n+1}$ from $x$ to $y$ such that $\widehat{S_n}\subseteq \widehat{S_{n+1}}\subseteq C_\infty$.

    Since $f$ is uniformly continuous, we can pick a real number $\delta_n$ such that:
    \begin{itemize}
        \item 
    \begin{equation}\label{cond_delta_n}
        0<\delta_n<\min\left\{\frac{\varepsilon_n}{6},\min\left\{\frac{d(x_i,x_j)}{2}\,\Bigg\vert\,i,j=0,\dots,m,\text{ and }x_i\neq x_j\right\}\right\};
    \end{equation}
    \item 
    
    for every $u,v\in X$, if $d(u,v)<\delta_n$, then $d(f(u),f(v))<\varepsilon_n/6$.
    \end{itemize}
    Since $\{C_k\}_{k=1}^\infty$ is a complete sequence of chains and by \eqref{hausconv}, there exists a natural number $k_n\in\mathbb{N}$  such that the chain
    $$C_{k_n}:\,x_0^{(k_n)},\,x^{(k_n)}_1,\dots,\,x_{m_{k_n}}^{(k_n)}$$
    is an $\varepsilon_n/6-$chain. We can take $k_n$ so large that $d_\mathcal{H}(\widehat{C_{k_n}},C_\infty)<\delta_n/2$. 
    Thus, since
    $$C_\infty\subseteq \overline{B}_{\delta_n/2}(\widehat{C_{k_n}}),$$
    every point $x_0,\dots,\,x_m\in\widehat{S_n}$ has a ``close enough'' point in $\widehat{C_{k_n}}$. More precisely,  
    for every $i=0,\dots,m$, since $\widehat{C_{k_n}}$ is finite, there is a point $x_{j_i}^{(k_n)}\in \widehat{C_{k_n}}$ such that
    \begin{equation}\label{def_assignment}
        d(x_{j_i}^{(k_n)},x_i)= \min_{z\in \widehat{C_{k_n}}}d(z,x_i),
    \end{equation}
    where $j_i\in\{0,\dots,m_{k_n}\}$ (note that it can be $x_{j_i}=x_i$).
    Moreover, by the definition of Hausdorff metric, we also have 
    \begin{equation}\label{dist-hau-proj}
        d(x_{j_i}^{(k_n)},x_i)\leq\delta_n/2<\delta_n.
    \end{equation}
    Note that, by \eqref{def_assignment}, if $x_i\in\widehat{C_{k_n}}$, then $x_{j_i}^{(k_n)}=x_i$; for instance $x_{j_0}^{(k_n)}=x_0=x$ and $x_{j_m}^{(k_n)}=x_m=y$.
    
    Now we want to make sure that the assignment
    \begin{equation}\label{assignment}
        \widehat{S_n}\ni x_i\mapsto x_{j_i}^{(k_n)}\in\widehat{C_{k_n}}
    \end{equation}
    is injective for $i\in\{0,\dots,m-1\}$. For this, recalling \eqref{cond_delta_n},    
    for every $i,i'\in\{0,\dots,m-1\}$ and $i'\neq i$, we have $d(x_{i'},x_{j_i}^{(k_n)})>\delta_n/2$, and therefore inequality \eqref{dist-hau-proj} implies injectivity of \eqref{assignment}.
    
    As a next step, we want to find points in $C_\infty$ that are ``close enough" to every point of $\widehat{C_{k_n}}$ that has not been assigned through \eqref{assignment}. For this observe that, by Hausdorff convergence, for every $$h\in\{1,\dots,m_{k_n}-1\}\setminus\{j_0,\dots,j_m\},$$ there exists $z_h\in C_\infty$ such that
    \begin{equation}
        d(z_h,x_h^{(k_n)})\leq\delta_n/2<\delta_n.
        \tag{\ref{dist-hau-proj}'}
        \label{dist-hau-proj2}
    \end{equation}

    In general, the indices $\{j_i\, ,\, i=0,\dots, m\}$ are not ordered by $i$, in the sense that it is not guaranteed that $j_0<j_1<\dots<j_m$. Let us thus define indices $\{i_p\, ,\, p\in 1,\dots, m-1\}$ so that $j_{i_1}<j_{i_2}<\dots<j_{i_{m-1}}$.
    
    Consider now the chain, supported on $C_\infty$, obtained by inserting between the points of $S_n$ (suitably reordered, with respect to their order on $S_n$, by the index $i_p$) the points $z_h$ that, by construction, belong to $C_\infty$ and are one by one suitably close to the unassigned points of $\widehat{C_{k_n}}$. More precisely, let us consider the points:
    \begin{multline*}
        x_0=x,\,z_1,\dots,\,z_{j_{i_1}-1},\,x_{i_1},\,z_{j_{i_1}+1},\dots,\,z_{j_{i_2}-1},\,x_{i_2},\,z_{j_{i_2}+1},\dots\\\dots,\,z_{j_{i_{m-1}}-1},\,x_{i_{m-1}},\,z_{j_{i_{m-1}}+1},\dots,\,z_{j_{i_m}-1},\,x_m=y.
    \end{multline*}
    
    Let us rename them as
    $$x_0',\,x_1',\dots,\,x_{m_{k_n}-1}',\,x_{m_{k_n}}'.$$

    Conditions (\ref{dist-hau-proj}) and (\ref{dist-hau-proj2}) yield $d(x_i^{(k_n)},x_i')<\delta_n$, for every $i\in\{0,\dots,m_{k_n}\}$, and thus we have also $d(f(x_i^{(k_n)}),f(x_i'))<\frac{\varepsilon_n}{3}$, for every $i\in\{0,\dots,m_{k_n}\}$.
    Therefore, for every $i=0,\dots,m_{k_n}-1$, we have:
    \begin{align*}
        d(f(x_i'), x_{i+1}') &\le d(f(x_i'), f(x_i^{(k_n)})) + d(f(x_i^{(k_n)}), x_{i+1}^{(k_n)}) + d(x_{i+1}^{(k_n)}, x_{i+1}') \\
        &< \frac{\varepsilon_n}{6} + \frac{\varepsilon_n}{6} + \delta_n \\
        &< \frac{\varepsilon_n}{2},
    \end{align*}
    where we use the fact that $C_{k_n}$ is an $\frac{\varepsilon_n}{6}$-chain and so $d(f(x_i^{(k_n)}), x_{i+1}^{(k_n)})<\frac{\varepsilon_n}{6}$.
    Thus,
    $$S_{n+1}:x_0',\,x_1',\dots,\,x_{m_{k_n}-1}',\,x_{m_{k_n}}'$$
    is an $\varepsilon_{n+1}$-chain, and, by construction, we also have
    $\widehat{S_n}\subseteq \widehat{S_{n+1}}$.

    \begin{rem}\label{inj__}
           Note that $\operatorname{card}(\widehat{S_{n+1}})\le m_{k_n}$ but in general  $\operatorname{card}(\widehat{S_{n+1}})\ne m_{k_n}$, because the assignment
    $$\widehat{C_{k_n}}\ni x_h^{(k_n)}\mapsto x_h'\in C_\infty$$
    is not guaranteed to be injective: two points of $\widehat{C_{k_n}}$ may share their closest point $z_{\bar{h}}$ in $C_\infty$ (notice that this problem has no easy fix if $z_{\bar{h}}$ is isolated in $C_\infty$).
    \end{rem}
    
    \begin{defi}\label{hauspro}
        We say that the collection of nested chains $\{S_n\}_{n=1}^\infty$ obtained above starting from the Hausdorff-converging complete sequence of chains $\{C_k\}_{k=1}^\infty$ is a \emph{Hausdorff projection} of $\{C_k\}_{k=1}^\infty$ over its limit, and indicate it by $\mathcal{H}(\{C_k\}_{k=1}^\infty)$.
    \end{defi}

\begin{rem}\label{S_n->C_infty}
    Note that, by construction, the chains
    $$S_{n+1}:\,x_0',\,x_1',\dots,\,x_{m_{k_n}}'$$
    and
    $$C_{k_n}:\,x_0^{(k_n)},\,x^{(k_n)}_1,\dots,\,x_{m_{k_n}}^{(k_n)}$$
    are such that, for every $i\in\{0,\dots,m_{k_n}\}$,
    $$d(x^{(k_n)}_i,x_i')\le\frac{\delta_n}{2}.$$
    Therefore, we have $d_\mathcal{H}(\widehat{S_{n+1}},\widehat{C_{k_n}})\le\frac{\delta_n}{2}$, hence
    $$d_\mathcal{H}(\widehat{S_{n+1}},C_\infty)\le d_\mathcal{H}(\widehat{S_{n+1}},\widehat{C_{k_n}})+d_\mathcal{H}(\widehat{C_{k_n}},C_\infty)\le \frac{\delta_n}{2}+\frac{\varepsilon_n}{6}<2\frac{\varepsilon_n}{6}<\varepsilon_n,$$
    and thus $\widehat{S_n}\xrightarrow{n\to\infty}C_\infty$, too. By Eq.\eqref{nest_haus}, this implies that $C_\infty=\overline{\cup_n \widehat{S_n}}$
\end{rem}

\begin{thm}\label{C_infty subsystem}
    Let $x,y\in X$ be such that $x\,\mathcal{C}\,y$ and let $\{C_n\}_{n=1}^\infty$ be a complete sequence of chains from $x$ to $y$ for the sequence of positive reals approaching 0 monotonically $\{\varepsilon_n\}_{n=1}^\infty$. After passing if needed to a convergent subsequence, let $C_\infty$ be the Hausdorff limit of $\{\widehat{C_n}\}_{n=1}^\infty$. Then we have: $$f(C_\infty)\cup\{x\}=C_\infty\cup\{f(y)\}.$$
\end{thm}
\begin{proof}
   Assume that $\widehat{C_n}\xrightarrow{n\to\infty}C_\infty$.
    First of all, let us prove $f(C_\infty)\cup\{x\}\subseteq C_\infty\cup\{f(y)\}$. Since $x\in\widehat{C_n}$ for every $n\in\mathbb N$, we have $x\in C_\infty$. So, take $z\in C_\infty\setminus\{y\}$ and let
    $$C_n:\,x_0^{(n)},\,x_1^{(n)},\dots,\,x_{m_n}^{(n)}.$$
    Then, there exists a sequence of points $\{x^{(n)}_{i_n}\}_{n=1}^\infty$, where $i_n<m_n$ eventually, such that $$x^{(n)}_{i_n}\xrightarrow{n\to\infty}z.$$
    Consider now the sequence of points $\{x^{(n)}_{i_n+1}\}_{n=1}^\infty$. Since $C_n$ is a $\varepsilon_n$-chain, we get $$d(f(x^{(n)}_{i_n}),x^{(n)}_{i_n+1})<\varepsilon_n.$$ As $f(x^{(n)}_{i_n})\xrightarrow{n\to\infty}f(z)$ by continuity of $f$, we have $x^{(n)}_{i_n+1}\xrightarrow{n\to\infty}f(z)$, too. Thus, $f(z)\in C_\infty$.

    In order to prove the other inclusion, let $z\in C_\infty\setminus\{x\}$. Hence, there is a sequence of points $\{x^{(n)}_{i_n}\}_{n=1}^\infty$ such that $x^{(n)}_{i_n}\xrightarrow{n\to\infty}z$, where $i_n>0$ eventually. The sequence $\{x^{(n)}_{i_n-1}\}_{n=1}^\infty$ admits at least one limit point, because $X$ is compact. Let $w\in X$ be one of the limit points; that is, there is a subsequence $\{x^{(n_k)}_{i_{n_k}-1}\}_{k=1}^\infty$ such that $x^{(n_k)}_{i_{n_k}-1}\xrightarrow{k\to\infty}w$. Since $d(f(x^{(n)}_{i_n-1}),x^{(n)}_{i_n})<\varepsilon_n$, then $f(x^{(n_k)}_{i_{n_k}-1})\xrightarrow{k\to\infty}z$. Therefore, every limit point of the sequence $\{x^{(n)}_{i_n-1}\}_{n=1}^\infty$ is a pre-image of $z$. Moreover, since $\widehat{C_n}\xrightarrow{n\to\infty}C_\infty$, we have $w\in C_\infty$.
\end{proof}
\begin{rem}
    From the previous result, it follows that if $f(y)\in C_\infty$, then $(C_\infty,f_\vert)$ is a subsystem of $(X,f)$; otherwise, if $f(y)\notin C_\infty$, it is still true that $(C_\infty\cup\overline{\mathcal{O}(y)},f_\vert)$ is a subsystem of $(X,f)$. Moreover, if $x\in f(C_\infty)$, then the system $(C_\infty\cup\overline{\mathcal{O}(y)},f_{\vert})$ is surjective.
\end{rem}

We are now ready to prove that we can always build nested chains in compact systems.

\begin{thm}\label{C=C_nested}
    Let $(X,f)$ be a compact topological dynamical system.
    Then $\mathcal{C}=\mathcal{C}_\subseteq$.
\end{thm}
\begin{proof}
    Let $x,y\in X$ be such that $x\,\mathcal{C}\, y$ and let $\{\delta_n\}_{n=1}^\infty$ be a sequence of positive real numbers converging to zero monotonically. Let $\{C_n\}_{n=1}^\infty$ be a complete sequence of chains from $x$ to $y$ for the sequence $\{\delta_n\}_{n=1}^\infty$. 
    Up to considering a sub-sequence, we can assume that $\{\widehat{C_n}\}_{n=1}^\infty$ converges in the Hausdorff metric to a closed set $C_\infty$. Then, the Hausdorff projection $\{S_n\}_{n=1}^\infty=\mathcal{H}(\{C_n\}_{n=1}^\infty)$ is a family of nested chains for the sequence $\left\{\frac{4\,\cdot\,\operatorname{diam}(X)}{2^n}\right\}_{n=1}^\infty$.
\end{proof}

\begin{rem}
    As Example \ref{non-compatto} below shows, this result generally fails in non-compact dynamical systems. In fact, in a non-compact space $X$, a complete sequence of chains $\{C_n\}_{n=1}^\infty$ might admit no  sub-family $\{C_{n_k}\}_{k=1}^\infty$ such that $\{\widehat{C_{n_k}}\}_{k=1}^\infty$ is a convergent sequence of closed subsets in the Hausdorff metric; and this is just what happens when one considers in Example \ref{non-compatto} a complete sequence of chains from $x$ to $y$. 
\end{rem}

\begin{example}\label{non-compatto}
    In this example we show that, in the general non-compact case, the relation $\mathcal{C}_\subseteq$ can be a strictly smaller than $\mathcal{C}$. 
    
    Let $X\subset \mathbb{R}^2$ be the set consisting of the two half-lines with origins in $x=(0,0)$ and $y=(0,1)$  and parallel to the $x$-axis and of the points $z^k_h$, where, for every $k\in\mathbb N$ and for every $h\in\{1,\dots,k\}$, the point $z^k_h$ has coordinates $\left(k,\frac{h}{k+1}\right)$. Note that $X$ is a closed, non-compact subset of $\mathbb R^2$.

    \begin{figure}[H]
        \centering
        \begin{tikzpicture}[scale=1]
            \fill (0,2) circle (1.5pt) node[above left] {$y$};
            \fill (0,-2) circle (1.5pt) node[below left] {$x$};
    
            \draw[line width=1pt] (0,2) -- (8.5,2);
            \draw[line width=1pt] (0,-2) -- (8.5,-2);

            \draw[line width=1pt] (9.5,2) -- (11,2);
            \draw[line width=1pt] (9.5,-2) -- (11,-2);

            \foreach \x in {8.8,9,9.2}
                \fill (\x,2) circle (1pt);

            \foreach \x in {8.8,9,9.2}
                \fill (\x,0) circle (1pt);

            \foreach \x in {8.8,9,9.2}
                \fill (\x,-2) circle (1pt);
    
            \foreach \x in {11.2,11.4,11.6}
                \fill (\x,2) circle (1pt);
            \foreach \x in {11.2,11.4,11.6}
                \fill (\x,0) circle (1pt);
            \foreach \x in {11.2,11.4,11.6}
                \fill (\x,-2) circle (1pt);
    
            \fill (1.5,0) circle (1.5pt) node[left]{$z^1_1$};
    
            \fill (3,0.66) circle (1.5pt) node[left]{$z^2_2$};
            \fill (3,-0.66) circle (1.5pt) node[left]{$z^2_1$};
    
            \fill (4.5,1) circle (1.5pt) node[left]{$z^3_3$};
            \fill (4.5,0) circle (1.5pt) node[left]{$z^3_2$};
            \fill (4.5,-1) circle (1.5pt) node[left]{$z^3_1$};
    
            \fill (6,1.2) circle (1.5pt) node[left]{$z^4_4$};
            \fill (6,0.4) circle (1.5pt) node[left]{$z^4_3$};
            \fill (6,-0.4) circle (1.5pt) node[left]{$z^4_2$};
            \fill (6,-1.2) circle (1.5pt) node[left]{$z^4_1$};
    
            \fill (7.5,1.32) circle (1.5pt) node[left]{$z^5_5$};
            \fill (7.5,0.66) circle (1.5pt) node[left]{$z^5_4$};
            \fill (7.5,0) circle (1.5pt) node[left]{$z^5_3$};
            \fill (7.5,-0.66) circle (1.5pt) node[left]{$z^5_2$};
            \fill (7.5,-1.32) circle (1.5pt) node[left]{$z^5_1$};

            \foreach \y in {1.8,1.6,1.4,1.2,1,0.8,0.6,0.4,0.2,0,-0.2,-0.4,-0.6,-0.8,-1,-1.2,-1.4,-1.6,-1.8}
                \fill (10.25,\y) circle (1.5pt);
    
            \node at (6,-3) {The space $X$ described in Example \ref{non-compatto}};
        \end{tikzpicture}
    \end{figure}

    Take $f$ to be the identity map $\operatorname{id}_X:X\longrightarrow X$ and consider the dynamical system $(X,f)$. It is not difficult to see that $x\,\mathcal{C}\,y$, while $x\hspace{-.3em}\not\hspace{-0.3em}{\mathcal{C}_\subseteq}\,y$. Indeed, given $\varepsilon>0$ and an $\varepsilon$-chain $C$ between $x$ and $y$, there exists $\delta>0$ sufficiently small such that there is not a $\delta$-chain $C'$ from $x$ to $y$ with $\widehat{C}\subseteq\widehat{C'}$. Since $\widehat{C}$ is a finite set, we can take $K$ to be the maximal natural number such that there exists a point $z^K_h\in C$ for some $h\in\{1,\dots,K\}$. Then, for every $0<\delta<\frac{1}{K+1}$, we have $B_\delta(z^K_h)\cap X=\{z^K_h\}$ and, therefore, for every $\delta$-chain $C'$ from $x$ to $y$, we have $z^K_h\notin\widehat{C'}$.
\end{example}


The dynamical relevance of an $\varepsilon$-chain is linked to its interpretation as a means of moving from one point to another using only corrections of size $<\varepsilon$. From this perspective, following a cyclic sub-path while going from $x$ to $y$ appears redundant, and its dynamical significance is therefore problematic. This informal observation has a technical counterpart, as in case of cyclic sub-chains the order induced by the index is not a linear order.
We will address the problem of obtaining a sequence of nested $\varepsilon_n$-chains in which each chain is acyclic. In this way, it will be natural to associate to every $\varepsilon$-chain a linear order dictated by the indexes along the chain. 
This will be done in the next section. 

\section{\texorpdfstring{$\mathcal{C}=\mathcal{C}_\preceq$}{C = C-preceq}}\label{sec5}

In this section we prove that $\mathcal{C}_\preceq$ always coincides with $\mathcal{C}$. 

To achieve our goal, we refine the construction of nested chains by means of a
transfinite pruning procedure, by which we build a family of closed sets
$\{C^\lambda\}_{\lambda<\omega_1}$ and nested chains
$\{S_n^\lambda\}_{n\in\mathbb{N}}$ supported on them. Heuristically, $C^\lambda$ is the region
where cyclic behaviour still survives after $\lambda$ rounds of pruning.
At each successor stage we remove cycles inside $C^\lambda$ and project again
onto a Hausdorff limit $C^{\lambda+1}\subseteq C^\lambda$, while at limit
ordinals we intersect. This process stabilizes at some countable ordinal
$\lambda$ such that $C^\lambda=C^{\lambda+1}$. At that stage the Hausdorff projection
can no longer identify distinct chain points, so acyclicity is preserved and we
obtain nested, acyclic chains that witness $x\,\mathcal{C}_\preceq\,y$.

The following lemma, ensuring stabilization of closed, nonempty, nested sets
at some countable level, is folklore in the setting of second-countable
spaces. We include a proof for completeness.

\begin{lem}\label{seq.nested.clo.sets.}
    Let $X$ be a second-countable space and let $\{C^\beta\}_{\beta<\omega_1}$ be a family of closed sets such that $C^{\beta+1}\subseteq C^\beta$ for every $\beta<\omega_1$. Then, there exists $\lambda<\omega_1$ such that $C^{\lambda+1}=C^\lambda$.
\end{lem}
\begin{proof}
    Take $\{B_n\}_{n=1}^\infty$ to be a countable base for the topology. Towards a contradiction, let us assume $C^{\beta+1}\subsetneq C^\beta$ for every $\beta<\omega_1$. Thus, for every $\beta<\omega_1$, there exists $x_\beta\in C^\beta\setminus C^{\beta+1}$. Furthermore, since $x_\beta\notin C^{\beta+1}$ and $C^{\beta+1}$ is closed, for every $\beta<\omega_1$ there exists a natural number $n_\beta\in\mathbb{N}$ such that $x_\beta\in B_{n_\beta}$ and $B_{n_\beta}\cap C^{\beta+1}=\varnothing$. Now, for every $\beta_1,\beta_2<\omega_1$ such that $\beta_1\neq\beta_2$, we have $n_{\beta_1}\neq n_{\beta_2}$. In fact, if $\beta_1<\beta_2$, then $x_{\beta_2}\in C^{\beta_2}\subseteq C^{\beta_1+1}$ and, since $B_{n_{\beta_1}}\cap C^{\beta_1+1}=\varnothing$ and $x_{\beta_2}\in B_{n_{\beta_2}}$, we have $B_{n_{\beta_1}}\neq B_{n_{\beta_2}}$. Therefore, we have found a contradiction since $\{B_n\}_{n=1}^\infty$ is countable.
\end{proof}

\begin{thm}\label{mainth}
    Let $(X,f)$ be a compact topological dynamical system. Then $\mathcal{C}=\mathcal{C}_\preceq$.
\end{thm}
\begin{proof}
    Let $x,y\in X$ be such that $x\,\mathcal{C}\,y$ and let $\{C_n\}_{n=1}^\infty$ be complete sequence of chains. Up to passing to a sub-family, suppose $\widehat{C_n}\xrightarrow{n\to\infty}C_\infty$ for some closed set $C_\infty$. We need to prove that we can find a sequence of nested chains $\{S_n\}_{n=1}^\infty$ such that, for every $n\in\mathbb{N}$, $S_n$ does not contain a cyclic sub-chain. 

    We will now define a family of pairs $\{(C^\lambda, \{S_n^\lambda\}_{n=1}^\infty)\}_{\lambda<\omega_1}$ indexed by countable ordinals, where $C^\lambda$ is a closed set such that $C^\lambda\subseteq C^{\lambda-1}$ and $\{S^\lambda_n\}_{n=1}^\infty$ is a family of nested chains such that, for every $n\in\mathbb{N}$ and for every ordinal number $\lambda$, $S^\lambda_n$ is an $\varepsilon_n$-chain and $\widehat{S^\lambda_n}\subseteq C^\lambda$, with $\varepsilon_n=\frac{4\cdot\operatorname{diam(X)}}{2^n}$.

    For $\lambda=0$, let $C^0$ and $\{S^0_n\}_{n=1}^{\infty}$ be respectively the set $C_\infty$ and the Hausdorff projection $\mathcal{H}(\{C_n\}_{n=1}^{\infty})$ over $C^0$. Then,
    \begin{itemize} 
        \item if $\lambda$ is a successor ordinal, consider the sequence of nested chains $\{S^{\lambda-1}_{n}\}_{n=1}^{\infty}$. For every $n\in\mathbb{N}$, we can modify the $\varepsilon_n$-chain $S^{\lambda-1}_n$ obtaining another $\varepsilon_n$-chain 
        \begin{equation}\label{acyclic_chain}
            ({S^{\lambda-1}_n})'
        \end{equation}
       with no cyclic sub-chain and such that $\widehat{({S^{\lambda-1}_n})'}\subseteq\widehat{S^{\lambda-1}_n}$. Now, since $\{\widehat{({S^{\lambda-1}_n})'}\}_{n=1}^{\infty}$ is a sequence of closed sets, up to passing to a sub-sequence, we can assume that $\{\widehat{({S^{\lambda-1}_n})'}\}_{n=1}^\infty$ is a converging sequence in the Hausdorff metric. Take $C^\lambda$ to be the Hausdorff limit of this sequence. By the inductive hypothesis we have $\widehat{{(S^{\lambda-1}_n})'}\subseteq\widehat{S^{\lambda-1}_n}\subseteq C^{\lambda-1}$ for every $n\in\mathbb{N}$; thus, we also have $C^\lambda\subseteq C^{\lambda-1}$. Then, let $\{S^\lambda_n\}_{n=1}^{\infty}$ be the Hausdorff projection $\mathcal{H}(\{({S^{\lambda-1}_n})'\}_{n=1}^{\infty})$.
        \item If $\lambda<\omega_1$ is a limit ordinal, take
        $$C^\lambda=\bigcap_{\beta<\lambda}C^\beta.$$
        Since $X$ is a compact space and,  $C^\beta\neq\varnothing$ for every $\beta<\lambda$, we have $C^\lambda\neq\varnothing$. Furthermore, as $C^{\beta+1}\subseteq C^\beta$ for every $\beta<\lambda$ by inductive hypothesis, we also have $C^\beta\xrightarrow{\beta\to\lambda}C^\lambda$ in the Hausdorff metric.
        Now, we need to define a sequence of nested chains $\{S^\lambda_n\}_{n=1}^{\infty}$ such that $S^\lambda_n$ is an $\varepsilon_n$-chain and $\widehat{S^\lambda_n}\subseteq C^\lambda$. For every $n\in\mathbb{N}$, there exist $\beta_n<\lambda$ such that $d_\mathcal{H}(C^\lambda,C^{\beta_n})<1/n$. Moreover, as said in Remark \ref{S_n->C_infty}, for every $n\in\mathbb{N}$, there exists $m_n\in\mathbb{N}$ such that $d_\mathcal{H}(\widehat{S^{\beta_n}_{m_n}},C^{\beta_n})<1/n$. Thus, the sequence of chains $\{S^{\beta_n}_{m_n}\}_{n=1}^{\infty}$ is such that
        $$d_\mathcal{H}(\widehat{S^{\beta_n}_{m_n}},C^\lambda)\le d_\mathcal{H}(\widehat{S^{\beta_n}_{m_n}},C^{\beta_n})+d_\mathcal{H}(C^\lambda,C^{\beta_n})<2/n,$$
        and so $\widehat{S^{\beta_n}_{m_n}}\xrightarrow{n\to\infty}C^\lambda$ in the Hausdorff metric. Then, we can take $\{S^\lambda_n\}_{n=1}^{\infty}$ to be the Hausdorff projection $\mathcal{H}(\{S^{\beta_n}_{m_n}\}_{n=1}^{\infty})$.
    \end{itemize}
    In this way we obtain a sequence $\{C^\lambda\}_{\lambda<\omega_1}$ of nested closed sets. 
    
    Since $X$ is a compact metric space, it is separable and second-countable, too. Then, by Lemma \ref{seq.nested.clo.sets.}, we know there is a countable ordinal $\lambda$ such that $C^\lambda=C^{\lambda+1}$. Therefore, the sequence of chains $\{({S^\lambda_n})'\}_{n=1}^{\infty}$ defined in \eqref{acyclic_chain} is such that $\widehat{({S^\lambda_n})'}\subseteq C^\lambda$ for every $n\in\mathbb{N}$, and moreover it also verifies
    \begin{equation}\label{hausconv2}        \widehat{({S^\lambda_n})'}\xrightarrow{n\to\infty}C^\lambda.
    \end{equation}

   So, by Eq.\eqref{hausconv2}, the Hausdorff projection at ordinal level $\lambda$ is taken onto a Hausdorff limit that already contains the projected chains. This ensures, as we will see, that the projected chains are acyclic provided the projecting ones are, because the “closest” point in the limit will in fact coincide with the projected point, which guarantees uniqueness.
    
    Indeed, let the chain $({S^\lambda_k})'$ consist of the points $$x_0^{(k)},\,x_1^{(k)},\dots,\,x^{(k)}_{m_k}.$$
    We now show that we can take a Hausdorff projection $\{S_n\}_{n=1}^\infty=\mathcal{H}(\{({S^\lambda_n})'\}_{n=1}^\infty)$ such that, for every $n$, we have that $S_n$ is an $\varepsilon_n$-chain from $x$ to $y$ with no cyclic sub-chain. We proceed inductively. Let $$S_1:\{0,1\}\longrightarrow C^\lambda$$ be the $\varepsilon_1$-chain consisting of just the points $x_0=x$ and $x_1=y$. Assume now that we have an $\varepsilon_n$-chain $S_n:x_0=x,\,x_1,\dots,\,x_m=y$ with no cyclic sub-chain such that $\widehat{S_n}\subseteq C^\lambda$. We prove that we can define an $\varepsilon_{n+1}$-chain $S_{n+1}$ such that $\widehat{S_n}\subseteq\widehat{S_{n+1}}\subseteq C^\lambda$ and such that it does not contain any cyclic sub-chain.

    Since $f$ is uniformly continuous, we can pick a real number $\delta_n$ such that:
    \begin{itemize}
        \item \begin{equation}\label{cond_delta_n_2}
            0<\delta_n<\min\left\{\frac{\varepsilon_n}{6},\min\left\{\frac{d(x_{i},x_{j})}{2}\,\Bigg\vert\,i,j=0,\dots,m,\text{ and }x_i\neq x_j\right\}\right\};
        \end{equation}
        \item if $d(u,v)<\delta_n$, then $d(f(u),f(v))<\varepsilon_n/6$.
    \end{itemize}

    By \eqref{hausconv2} and since $\{({S_n^\lambda})'\}_{n=1}^{\infty}$ is a complete sequence of chains, there exists a natural number $k_n\in\mathbb N$ such that $d_\mathcal{H}(\widehat{({S_{k_n}^\lambda})'},C^\lambda)<\delta_n/2$ and ${(S_{k_n}^\lambda})'$ is an $\varepsilon_n/6$-chain.

    Thus, since $C^\lambda\subseteq \overline{B}_{\delta_n/2}(\widehat{({S^\lambda_{k_n}})'})$, for every $i=0,\dots,m$ there are points $x_{j_i}^{(k_n)}$ of the chain $({S^\lambda_{k_n}})'$ such that
    $$d(x_{j_i}^{(k_n)},x_i)=\min_{z\in \widehat{({S_{k_n}^\lambda})'}} d(z,x_i),$$
    
    where $j_i\in\{0,\dots,m_k\}$ (note that it can be $x_{j_i}=x_i$). Moreover, by definition of Hausdorff metric, we have also $$d(x_{j_i}^{(k_n)}, x_i)\leq\delta_n/2<\delta_n.$$
    Condition \eqref{cond_delta_n_2} implies that, for every $i'\in\{0,\dots,m\}$ and $i'\neq i$, we have $d(x_{i'},x_{j_i}^{(k_n)})>\delta/2$; and therefore the assignment
    \begin{equation}\label{assignement_preceq}
        \widehat{S_n}\ni x_i\mapsto x_{j_i}^{(k_n)}\in \widehat{({S^\lambda_{k_n}})'}
    \end{equation}
    is injective.

    As done when defining the Hausdorff projection, we now want to find points in $C^\lambda$ that are ``close enough" to every point of $\widehat{({S^\lambda_{k_n}})'}$ that has not been assigned through \eqref{assignement_preceq}. Remember, though, that $\widehat{({S^\lambda_{k_n}})'}$ is a subset of the Hausdorff limit $C^\lambda$, and thus we can define, for every $h\in\{1,\dots,m_{k_n}-1\}\setminus\{j_0,\dots,j_m\}$, the points $z_h$ simply as $z_h=x_h^{(k_n)}$.

    The assignment 
    \begin{equation}\label{assignment_preceq_2}
        \widehat{({S^\lambda_{k_n}})'}\ni x_h^{(k_n)}\mapsto z_h\in C^\lambda
    \end{equation}
    is thus injective for $h\notin\{j_0,\dots,j_m\}$.

    Since the indices $\{j_i\mid i=0,\dots, m\}$ are not necessarily ordered by $i$, in the sense that it is not guaranteed that $j_0<j_1<\dots<j_m$, we define indices $\{i_p\mid p\in 1,\dots, m-1\}$ such that $j_{i_1}<j_{i_2}<\dots<j_{i_{m-1}}$.
    
    Consider now the chain, supported on $C^\lambda$, obtained by inserting between the points of $S_n$ (suitably reordered, with respect to their order on $S_n$, by the index $i_p$) the points $z_h$ that, by construction, belong to $C^\lambda$ and are one by one suitably close to the unassigned points of $\widehat{({S_{k_n}^\lambda})'}$. More precisely, let us consider the chain
    \begin{multline*}
        S^*:=x_0=x,\,z_1,\dots,\,z_{j_{i_1}-1},\,x_{i_1},\,z_{j_{i_1}+1},\dots,\,z_{j_{i_2}-1},\,x_{i_2},\,z_{j_{i_2}+1},\dots\\\dots,\,z_{j_{i_{m-1}}-1},\,x_{i_{m-1}},\,z_{j_{i_{m-1}}+1},\dots,\,z_{j_{i_m}-1},\,x_m=y.
    \end{multline*}

We observe that $S^*$ does not contain any cyclic sub-chain. Indeed, since the matching map
$$
\widehat{S_n}\longrightarrow \widehat{({S^\lambda_{k_n}})'}
$$
defined in \eqref{assignement_preceq} is injective and the chain $({S^\lambda_{k_n}})'$ is  acyclic, every point in $\widehat{S_n}$ and every new point $z_h$ has a unique ``ancestor'' in $\widehat{({S^\lambda_{k_n}})'}$. This ensures that no repetition can be created when building $S^*$. More precisely:
    \begin{itemize}
        \item for every $h,l\in\{1,\dots,m_k-1\}\setminus\{j_0,\dots,j_m\}$ with $h\neq l$, we have $z_h\neq z_l$, because of the injectivity of \eqref{assignment_preceq_2};
        \item for every $h,l\in\{0,\dots,m\}$ with $h\neq l$, we have $x_h\neq x_l$, because $S_n$ is acyclic by inductive hypothesis;
        \item  for every $h\in\{1,\dots,m_k-1\}\setminus\{j_0,\dots,j_m\}$ and for every $l\in\{0,\dots,m\}$, we have $z_h\neq x_l$. In fact, in this case, $z_h$ is a point in ${(S^\lambda_{k_n}})'$ that was not assigned through \eqref{assignement_preceq}. If we had $z_h=x_l$, then $z_h$ would have been such that
        $$
        d(x_l,z_h)=0=\min_{z\in\widehat{({S^\lambda_{k_n}})'}}d(x_l,z),
        $$
        and then \eqref{assignement_preceq} would assign $z_h$ to $x_l$, a contradiction.
    \end{itemize}
    Finally, we now show that $S^*$ is an $\varepsilon_n/2$-chain. For simplicity, let us rename the elements of the chain $S^*$ as
    $$S^*:\,x_0',\,x_1',\dots,\,x_{m^*-1}',\,x_{m_{k_n}}'.$$

    By construction, we have $d(x_i^{(k_n)},x_i')<\delta_n$, for every $i\in\{0,\dots,m_{k_n}\}$, and thus we have also $d(f(x_i^{(k_n)}),f(x_i'))<\frac{\varepsilon_n}{3}$, for every $i\in\{0,\dots,m_{k_n}\}$.
    Therefore, for every $i=0,\dots,m_{k_n}-1$, we have:
    \begin{align*}
        d(f(x_i'), x_{i+1}') &\le d(f(x_i'), f(x_i^{(k_n)})) + d(f(x_i^{(k_n)}), x_{i+1}^{(k_n)}) + d(x_{i+1}^{(k_n)}, x_{i+1}') \\
        &< \frac{\varepsilon_n}{6} + \frac{\varepsilon_n}{6} + \delta_n \\
        &< \frac{\varepsilon_n}{2},
    \end{align*}
    where we use the fact that $({S_{k_n}^\lambda})'$ is an $\frac{\varepsilon_n}{6}$-chain and so $d(f(x_i^{(k_n)}), x_{i+1}^{(k_n)})<\frac{\varepsilon_n}{6}$.
    Thus, we can define the desired $\varepsilon_{n+1}$-chain as $ S_{n+1}:= S^*$, which concludes the proof.
    
\end{proof}

The following definition imposes  order-compatibility for nested, acyclic chains. Note that this has to be required explicitly, as it cannot be ensured just by passing to subsequences, as it was possible with weak order-compatibility. In general indeed, as recalled in Remark \ref{diag_}, a sequence of nested, linearly ordered finite sets has a weakly order-compatible subsequence, but does not have an order-compatible subsequence. 

\begin{defi}
   Let $(X,f)$ be a dynamical system. For $x,y\in X$ such that $x\,\mathcal{C_\preceq}\,y$, we say that the complete sequence of nested chains $\{C_n\}_{n=1}^\infty$ is a \emph{order-compatible} if each $C_n$ has no cyclic sub-chains and for every $n\in \mathbb{N}$, the linear order $\le_{n+1}$ on $C_{n+1}$ extends the linear order $\le_n$ on $C_n$. 
\end{defi}

\begin{rem}
    If $\{C_n\}_{n=1}^\infty$ is a complete sequence of order-compatible nested chains from $x$ to $y$, for some $x,y\in X$, then, for every $n\in\mathbb N$, the chain $C_{n+1}$ enriches the chain $C_n$ without changing the order in which the points of $C_n$ appear. Moreover, for every $z,w\in\cup_n\widehat{C_n}$, the sequence $\{C_n\}_{n=1}^\infty$ restricted between $z$ and $w$ is still a complete sequence of order-compatible nested chains.
    Let us make this more precise.

    Without lost of generality, suppose $z\le_nw$ for every $n\in\mathbb N$ such that $z,w\in\widehat{C_n}$. Let $N=\min\{n\in\mathbb N:z,w\in\widehat{C_n}\}$ and, for every $n\ge N$, let
    $$
    C_n: x_0^{(n)}=x,x^{(n)}_1,\dots,\,x^{(n)}_{i_n-1},\,x^{(n)}_{i_n}=z,\,x^{(n)}_{i_n+1},\dots,\,x^{(n)}_{j_n-1},\,x^{(n)}_{j_n}=w,\,x^{(n)}_{j_n+1},\dots,\,x^{(n)}_{m_n}=y.$$
    We want to prove that $\{C_n^{z\to w}\}_{n=N}^\infty$ is a sequence of nested chains, where 
    $$
    C_n^{z\to w}:\,x_{i_n}^{(n)},\,x_{i_n+1}^{(n)},\dots,\,x_{j_n-1}^{(n)},\,x_{j_n}^{(n)}.
    $$
    
    Towards a contradiction, suppose there exists an index $i_n<h_n<j_n$ and a natural number $K\in\mathbb N$ such that the point
    $$
    x^{(n)}_{h_n}\notin\widehat{C_{n+K}^{z\to w}}.
    $$
    Since $\widehat{C_n}\subseteq \widehat{C_{n+K}}$, the point $x_{h_n}$ must appear either before $z$ or after $w$ in the chain $C_{n+K}$, which means that
    $$z\le_n x_{h_n}^{(n)}\quad\text{and}\quad x_{h_n}^{(n)}\le_{n+K} z$$
    or
    $$x_{h_n}^{(n)}\le_n w\quad\text{and}\quad w\le_{n+K} x_{h_n}^{(n)}$$ but this contradicts the hypothesis that $\{C_n\}_{n=1}^\infty$ is a sequence of order-compatible nested chains from $x$ to $y$ because the linear order $\le_{n+K}$ on $C_{n+K}$ does not extend $\le_n$ on $C_n$.
\end{rem}

We conclude this section with two results that will prove useful in the following.

\begin{lem}\label{lem:ONC-implies-chain}
Let $\{S_n\}_{n=1}^{\infty}$ be a complete sequence of order-compatible nested chains from $x$ to $y$
and let $S:=\bigcup_n \widehat{S_n}$. Let $z,w\in S$ be
distinct. Suppose that, for all sufficiently large $n$, $z$ appears
before $w$ in $S_n$. Then $z\,\mathcal{C}\,w$.
\end{lem}

\begin{proof}
Let $\{\varepsilon_n\}_{n=1}^\infty$ be a sequence with $\varepsilon_n\to 0$ such
that each $S_n$ is an $\varepsilon_n$-chain from $x$ to $y$. If $z,w \notin \{x,y\}$, there exists $N$ such that for all $n\ge N$
both $z$ and $w$ appear in $S_n$ and $z$ precedes $w$. For each $n\ge N$
let
$$
   S_n^{z\to w} :z=x_{i_n}^{(n)},\ x_{i_n+1}^{(n)},\dots,x_{j_n}^{(n)}=w
$$
denote the sub-chain of $S_n$ from the occurrence of $z$ to the occurrence of $w$. Since $S_n$ is an $\varepsilon_n$-chain, so
is $S_n^{z\to w}$. As $\varepsilon_n\to 0$, this shows that for every
$\varepsilon>0$ there is an $\varepsilon$-chain from $z$
to $w$. Hence $z\,\mathcal{C}\,w$. The thesis trivially holds if one of two points coincides with $x$ or with $y$.
\end{proof}
\begin{lem}\label{successor}
Let $\{S_n\}_{n=1}^{\infty}$ be a complete sequence of order-compatible nested chains from $x$ to $y$
and let $S:=\bigcup_n \widehat{S_n}$. Let $z,w\in S\setminus\{x,y\}$ be
distinct. Suppose that $w$ is the successor of $z$ in $(S\setminus\{x,y\},\le_\infty)$, then $f(z)=w$.
\end{lem}
\begin{proof}
Let $N=\min\{n\in\mathbb{N}: z,w\in\widehat{S_n}\}$. We would like to prove that, for every $n\geq N$, the chain $S_n$ is given by $$S_n: x_0^{(n)}=x,x^{(n)}_1,\dots,x^{(n)}_{i_n}=z,x^{(n)}_{i_n+1}=w,\dots,x^{(n)}_{m_n}=y.$$
In fact, if there was $N'>N$ such that $z\le_{N'} a\le_{N'} w $, where $a$ is different from $z$ and $w$, since $\{S_n\}_{n=1}^\infty$ is a complete sequence of order-compatible nested chains, we would have that $z\le_n a\le_n z$ for every $n>N'$, and so $x\le_\infty a\le_\infty w$, which is in contradiction with the fact that $w$ is the successor of $z$ in $(S\setminus\{x,y\},\le_\infty)$.

Now, since $\{S_n\}_{n=1}^\infty$ is a complete sequence of chains, we have that $d(f(z),w)<\varepsilon$ for every $\varepsilon>0$ and then $f(z)=w$.
\end{proof}
\section{The emergent order spectrum}
\begin{defi}
    We indicate by:
    \begin{enumerate}
        \item [i)] $\omega$ the first infinite ordinal;
        \item [ii)] $\omega^*$ the reverse order-type of $\omega$, that is the order type of non-positive integers under the usual $\leq$ relation;
        \item[iii)] $\zeta=\omega^*+\omega$ the order type of the integer numbers $\mathbb Z$;
        \item[iv)] $\eta$ the countable dense order-type without extrema; 
        \item[v)] $\text{Lin}_{\aleph_0}$ the set of all countable linear order-types.
    \end{enumerate} 
\end{defi}

We will use standard ordinal addition and multiplication throughout; e.g.\ $\zeta\cdot\omega$ denotes the (countable) sum of $\omega$ many copies of $\zeta$.





\begin{defi}
    Let $(X,f)$ be a topological dynamical system. Consider two points $x,y\in X$. We let $\Omega_f(x,y)$ be the set of all the order-types $\beta$ such that there is a complete sequence of order-compatible nested chains $\{S_n\}_n$ from $x$ to $y$ such that, setting $S:=\cup_n \widehat{S_n}$, we have that $(S\setminus\{x,y\},\leq_\infty)$ has order type $\beta$. We will omit the subscript and write simply $\Omega(x,y)$ if the map is clear from the context. 
    We set $\Omega(x,y)=\varnothing$ if there is no complete sequence of order-compatible nested chains from $x$ to $y$, so in particular, if $(x,y)\notin\mathcal{C}$, then $\Omega(x,y)=\varnothing$.

    We will call $\Omega(x,y)$ \emph{the Emergent Order Spectrum} (EOS) of $(x,y)$.
    The map:
    \begin{equation}
        \Omega_f:X^2\ni (x,y)\mapsto \Omega_f(x,y)\in \mathcal{P}(\text{Lin}_{\aleph_0})
    \end{equation}
    will be called the \emph{EOS map} of the dynamical system $(X,f)$.

    Given a countable order-type $\xi$ and $x\in X$, we set $$[\xi](x):=\{y\in X\, |\, \xi\in\Omega(x,y)\}.$$

\end{defi}
\begin{rem}
    Notice that $\Omega(x,y)=\varnothing$ here means that there are no emergent orders between $x$ and $y$ (this happens always when $x$ is not chain-related to $y$), whereas $\Omega(x,y)\ni\varnothing$ indicates the fact that $f(x)=y$ so that we can define a complete sequence of order-compatible nested chains $\{S_n\}_{n=1}^{\infty}$ setting, for each $n\in\mathbb N$,
    $$S_n:\,x,\,y.$$
    In fact, in this case, we have $S\setminus\{x,y\}=\varnothing$.
\end{rem}

The following are elementary properties of the EOS. Note that 6. is a consequence of 4. and 5. The other ones are proven in the following.
\begin{enumerate}
\item (Metric independence) $\Omega$ is independent of the choice of compatible metric on $X$.
\item (Conjugacy invariance) If $h$ topologically conjugates $f$ to $g$, then
$\Omega_f(x,y)=\Omega_g(h(x),h(y))$.
\item (Orbit detection) A finite ordinal $k$ lies in $\Omega(x,y)$
iff $f^{k+1}(x)=y$.
\item (Periodicity detection) There are finite ordinals $k$ and $k'$ such that $k\in\Omega(x,y)$ and $k'\in\Omega(y,x)$ iff $x$ is a periodic point (if $x=y$ is a fixed point, then $k$ and $k'$ are the empty order).
\item (Pure recurrence detection) The first infinite ordinal $\omega$ is in $\Omega(x,y)$ iff $(x,y)\in \mathcal{R}\setminus\mathcal{O}$.
\item (Limit set detection) 
$$
y\in\omega(x)\quad\Longleftrightarrow\quad
\Big(\,\omega\in\Omega(x,y)\,\Big)\ \ \text{or}\ \ \Big(\exists\,m,n \in\mathbb{N}:\ m\in\Omega(x,y)\ \text{and}\ n\in\Omega(y,x)\Big).
$$
\item (Non-wandering relation witnessing) The pair $(x,y)$ belongs to $\mathcal{N}\setminus\mathcal{R}$ if the infinite ordinal $\zeta$ is in $\Omega(x,y)$.
\end{enumerate}

\begin{lem}\label{lem:subsequence}
    Let $(X,f)$ be a compact metrizable dynamical system, and let $x,y\in X$ be such that $x\,\mathcal{C}\,y$. Let $\{S_n\}_{n=1}^{\infty}$ be a complete sequence of order-compatible nested  $\varepsilon_n$-chains from $x$ to $y$ such that, setting $S:=\bigcup_n \widehat{S_n}$, the ordered limit set $(S\setminus\{x,y\},\leq_\infty)$ has order type $\beta$. Then, for every sub-family $\{S_{n_k}\}_{k=1}^{\infty}$ that is a complete sequence of order-compatible nested $\delta_k$-chains from $x$ to $y$, and setting $S'=\bigcup_k\widehat{S_{n_k}}$, we have that $S=S'$ and the ordered limit set $(S'\setminus\{x,y\},\le_\infty')$ has order type $\beta$, too.
\end{lem}
\begin{proof}
    We now prove $S=S'$. Trivially, we have $S'\subseteq S$. Take now $z\in S$. Thus, since $\{S_n\}_{n=1}^{\infty}$ is a complete sequence of nested chains, there exists $N\in\mathbb N$ such that $z\in\widehat{S_n}$, for every $n>N$. Since $\{S_{n_k}\}_{k=1}^{\infty}$ is a sub-sequence of $\{S_n\}_{n=1}^{\infty}$, there exists $K\in\mathbb N$ such that $n_K>N$. Hence, $z\in\widehat{S_{n_K}}\subseteq S'$       
    
    To prove the statement, we now show that, for every $z,w\in S\setminus\{x,y\}$ such that $z\le_\infty w$, we have $z\le_\infty'w$. In fact, $z\le_\infty w$ if and only if $z\leq_n w$ for every $n\in\mathbb N$ such that $z,w\in\widehat{S_n}$. Then, for every $k\in\mathbb{N}$ such that $z,w\in\widehat{S_{n_k}}$, we have $z\leq_{n_k} w$ in $(\widehat{S_{n_k}}\setminus\{x,y\},\leq_{n_k})$ and so $z\leq_\infty' w$.
\end{proof}

\begin{thm}\label{cor:sequence-indipendence}
    Let $(X,f)$ be a compact dynamical system. Then, for every $x,y\in X$, the ordering $\beta$ is in the set $\Omega(x,y)$ if, for every sequence of positive real numbers $\{\varepsilon_n\}_{n=1}^{\infty}$ that goes to zero monotonically, there exists a sequence of order-compatible nested $\varepsilon_n$-chains $\{S_n\}_{n=1}^{\infty}$ such that, setting $S=\bigcup_n\widehat{S_n}$, the ordered limit set $(S\setminus\{x,y\},\le_\infty)$ has order type $\beta$.
\end{thm}

\begin{proof}
    Take $x,y\in X$. If $\Omega(x,y)=\varnothing$ there is nothing to prove.
    Let thus $x,y\in X$ be such that $x\,\mathcal{C}_\preceq\,y$ and let moreover be $\beta\in\Omega(x,y)$. Thus, there exist a sequence of positive real numbers $\{\delta_n\}_{n=1}^\infty$ that tends to zero monotonically and a sequence of order-compatible nested chains $\{S_n\}_{n=1}^\infty$ from $x$ to $y$ for the sequence $\{\delta_n\}_{n=1}^\infty$ such that the ordered limit set $(S\setminus\{x,y\},\le_\infty)$ has order type $\beta$, where $S=\bigcup_n\widehat{S_n}$.

    Take now another sequence of positive real numbers $\{\varepsilon_k\}_{k=1}^{\infty}$ that tends to zero monotonically. By Remark \ref{canon_}, we know that we can extract a sub-sequence $\{S_{n_k}\}_{k=1}^{\infty}$ such that $\{S_{n_k}\}_{k=1}^{\infty}$ is a sequence of order-compatible nested chains from $x$ to $y$ for the sequence $\{\varepsilon_k\}_{k=1}^{\infty}$. By Lemma \ref{lem:subsequence} we also have that the order type of the ordered limit set $$\left(\left(\bigcup_k\widehat{S_{n_k}}\right)\setminus\{x,y\},\le_\infty\right)$$ is $\beta$, so we are done.
\end{proof}

\begin{rem}\label{indep_eps_}
    By theorem \ref{cor:sequence-indipendence}, given two points $x,y\in X$ in a compact dynamical system $(X,f)$, the set $\Omega(x,y)$ is independent of the sequences $\{\varepsilon_n\}_{n=1}^{\infty}$ we chose to construct the sequence of order-compatible nested chains $\{S_n\}_{n=1}^{\infty}$ to obtain a limit set $$\left(\left(\bigcup_n\widehat{S_{n}}\right)\setminus\{x,y\},\le_\infty\right)$$ with a specific order type.   
\end{rem}



$$
$$



Let us now explore the interplay between the EOS of a pair of points and the kind of recurrence that may occur between them. 
Let us start by the invariance result.

\begin{thm}\label{conj__}
The EOS map $\Omega$ is invariant under topological conjugacy. 

Precisely, if $(X,f)$ and $(Y,g)$ are topologically conjugate through the homeomorphism $h:X\to Y$, then, setting $$H=h\times h:X^2\ni (x,y)\mapsto (h(x),h(y))\in Y^2,$$ we have $$\Omega_f=\Omega_g\circ\  H\text{ \emph{and} }\ \Omega_g= \Omega_f\circ H^{-1}.$$
\end{thm}
\begin{proof}
Let $x,y\in X$ and let a countable order-type $\beta$\ be such that $\beta\in \Omega_f (x\, ,y)$. We would like to prove that $\beta\in\Omega_g(h(x),h(y))$. In fact, consider $\{\delta_n\}_{n=1}^\infty$  a sequence of positive real numbers converging to zero and, for every $n\in\mathbb{N}$, since $h$ is a uniformly continuous function, there exists $\varepsilon_n>0$ such that, for every $z_1\, ,z_2 \in X$, $$d_X(z_1\,,\,z_2)<\varepsilon_n \Rightarrow\, d_Y (h(z_1)\, ,\, h(z_2))<\delta_n.$$ Since it is always possible to take $\varepsilon_n<\delta_n$ for every $n\in \mathbb{N}$, also $\{\varepsilon_n\}_{n=1}^\infty$ is  a sequence of positive real numbers converging to zero. Let $\{S_n\}_n$ be a sequence of order-compatible nested chains from $x$ to $y$ for the sequence $\{\varepsilon_n\}_{n=1}^\infty$ such that, setting $S:=\cup_n \widehat{S_n}$, the ordered limit set $(S\setminus\{x,y\},\leq_\infty)$ has order type $\beta$. If $S_n: x_0,\,x_1,\dots,\,x_m$, setting $D_n=h\circ S_n$, we have trivially $\widehat{D_n}\subseteq \widehat{D_{n+1}}$. In addition, $\{D_n\}_{n=1}^{\infty}$ is a sequence of order-compatible nested chains from $x$ to $y$ for the sequence $\{\delta_n\}_{n=1}^\infty$. In fact, we have $D_n:\,h({x_0}),h({x_1}),\dots,h({x_m})$ and, for every $i=0,\dots,m-1$,  $$d_Y(g(h({x_i})),h(x_{i+1}))=d_Y(h(f({x_i})),h(x_{i+1}))<\delta_n$$ where in the last inequality we used the fact that $d_X(f({x_i}),x_{i+1})<\varepsilon_n.$ Then, the limit ordered set $(D\setminus\{h(x),h(y)\},\leq_\infty)$, where $D=\bigcup_n \widehat{D_n}$, has order type $\beta$. Finally, performing the same construction and using the fact that $h$ is a homeomorphism, one proves that $\Omega_g= \Omega_f\circ H^{-1}$
\end{proof}

\begin{rem}\label{rem:chain-indipendence-metric}
From the previous result, it follows immediately (using the identity as the conjugating homeomorphism) that the EOS is metric-invariant. More precisely, let $(X,f)$ be a compact metrizable dynamical system, and let $d$ and $d'$ be two compatible metrics on $X$. Then for every $x,y\in X$,
$$
\Omega_{d}(x,y)\;=\;\Omega_{d'}(x,y).
$$
In particular, $\Omega$ does not depend on the chosen compatible metric.
\end{rem}

\begin{thm}\label{finite-ordinal}
    Let $x,y\in X$. Then 
    $k\in\Omega(x,y)$ (for $k\in\mathbb{N}_0$) if and only if $x\,\mathcal{O}\,y$.
\end{thm}
\begin{proof}

    Let $k\in\Omega(x,y)$, and let $\{S_n\}_{n=1}^{\infty}$ be a sequence of order-compatible nested chains from $x$ to $y$ for a sequence $\{\varepsilon_n\}_{n=1}^\infty$ such that, setting $S=\bigcup_n\widehat{S_n}$, the ordered limit set $(S\setminus\{x,y\},\le_\infty)$ has order type $k$. Then, we can enumerate the elements of $S\setminus\{x,y\}$ as $x_1,\,x_2,\dots,\,x_k$, where $x_i\le_\infty x_{i+1}$ for every $i\in\{1,\dots,k-1\}$.
    Since $S\setminus\{x,y\}$ is finite, by the definition of $\le_\infty$,  there exists $N\in\mathbb{N}$ such that $$S_n:x_0=x,\,x_1,\dots,x_k,\,x_{k+1}=y $$
    for every $n>N$ and this can be true only if $f(x_i)=x_{i+1}$ for every $i=0,\dots,k$, which means that $f^{k+1}(x)=y$.

    Vice versa, if $x\,\mathcal{O}\,y$, let $k\in\mathbb{N}$ be the minimal natural number such that $f^{k+1}(x)=y$. Then, $x,\,f(x),\dots,\,f^{k+1}(x)=y$ is an $\varepsilon$-chain for every $\varepsilon>0$ which does not contain any cyclic sub-chain. Hence, taking for every $n\in\mathbb N$
    $$S_n:\,x,\,f(x),\dots,\,f^{k+1}(x)=y,$$ we obtain a linearly ordered limit set $(S\setminus\{x,y\},\le_\infty)$ isomorphic to the ordinal $k$, where $S=\bigcup_n\widehat{S_n}$. Thus, $k\in\Omega(x,y)$.
\end{proof}

\begin{rem}
    Note that, if there exist two distinct finite ordinals $k$ and $k'$ such that $k\in\Omega(x,y)$ and $k'\in\Omega(y,x)$, then $x$ is a periodic point with period $k+k'+2$.
\end{rem}
\begin{thm}\label{Omega_primo_R}
    For $x,y\in X$, the first infinite ordinal $\omega$ is in $\Omega(x,y)$ if and only if $(x,y)\in\mathcal{R}\setminus\mathcal{O}$.
\end{thm}
\begin{proof}
    Suppose $(x,y)\in\mathcal{R}\setminus\mathcal{O}$. Given a sequence of positive real numbers $\{\varepsilon_n\}_{n=1}^{\infty}$ converging to 0 monotonically,
    we can consider a strictly increasing sequence of natural numbers $\{k_n\}_{n=1}^{\infty}\subseteq\mathbb
    N$ such that $d(f^{k_n}(x),y)<\varepsilon_n$. Set
    $$S_n:\,x,\,f(x),\dots,\,f^{k_n-1}(x),\,y.$$ 
    It is not difficult to see that $\{S_n\}_{n=1}^{\infty}$ is a sequence of order-compatible nested chains, and the order-type of the ordered limit set $(S\setminus\{x,y\},\le_\infty)$ (where $S=\bigcup_n\widehat{S_n}$) is $\omega$.

 Suppose now that $\omega\in\Omega(x,y)$ and let $\{S_n\}_{n=1}^{\infty}$ be a complete sequence of order-compatible nested chains such that the ordered limit set $(S\setminus\{x,y\},\le_\infty)$ has order type $\omega$, where $S=\bigcup_n\widehat{S_n}$.
     Notice that $(S\setminus\{x,y\},\le_\infty)$ admits a minimal element $z$. This can be true only if, ultimately, the second element in the chains $S_n$ is $z$; which means that
     $$d(f(x),z)<\varepsilon$$
     for every $\varepsilon>0$, and thus, $z=f(x).$
     Since $(S\setminus\{x,y\},\le_\infty)$ is isomorphic to $\omega$, for every element in $S\setminus\{x,y\}$, there is a successor. In particular, note that, if $w$ is the successor of $v$, then, by Lemma \ref{successor}, we have $f(v)=w$.
     Therefore, we have
     $$S\setminus\{x,y\}=\mathcal{O}(z)\cup\{z\}=\mathcal{O}(x).$$
    Since $\omega\in\Omega(x,y)$, we have $y\notin\mathcal{O}(x)$. In fact, otherwise, we would have $\operatorname{card}(S)<\infty$. Since $\{S_n\}_{n=1}^{\infty}$ is a complete sequence of chains from $x$ to $y$, we have
    $$y\in\overline{\mathcal{O}(x)},$$
    which means $x\,\mathcal{R}\,y$. 
\end{proof}




\begin{thm}\label{zetaN}
    Let $(X,f)$ be a dynamical system where $f:X\longrightarrow X$ is a homeomorphism and let $x,y\in X$ be such that  $\zeta\in\Omega(x,y)$. Then $(x,y)\in\mathcal{N}$.
\end{thm}

\begin{proof}
    Let $\{S_n\}_{n=1}^{\infty}$ be a complete sequence of order-compatible nested chains from $x$ to $y$ such that the ordered limit set $(S\setminus\{x,y\},\le_\infty)$, where $S=\bigcup_n\widehat{S_n}$, has order-type $\zeta$. Then, for every element in $S\setminus\{x,\,y\}$, there are a successor and a predecessor in $(S\setminus\{x,y\},\le_\infty)$. Applying Lemma \ref{successor}, we have that if $w$ is the successor of $z$ in $S\setminus\{x,y\}$ then $f(z)=w$. Therefore, for every $w\in S\setminus\{x,y\}$, there exists $z\in S\setminus\{x,y\}$ such that $f(z)=w$ and, moreover, $f(w)\in S$, too.
    This assures that there exists $\{z_n\}_{n=1}^\infty\subseteq S$
    such that $y\in \overline{\mathcal{O}(z_1)}$, $z_{n+1}\,\mathcal{O}\, z_n$ and $d(f(x),z_n)<\frac{1}{n}$.
    By the fact that $f$ is a homeomorphism, the sequence $\{f^{-1}(z_n)\}_n$ is such that, for every $n\in\mathbb N$,
    $$d(x,f^{-1}(z_n))\xrightarrow{n\to\infty}0\quad,\quad f^{-1}(z_{n+1})\,\mathcal{O}\,f^{-1}(z_n)\quad\text{and}\quad y\in\overline{\mathcal{O}(f^{-1}(z_n))}.$$
    Thus, we have $(x,y)\in\mathcal{N}$.
\end{proof}
The converse of the previous theorem is false, in general (see Theorem \ref{denjo_}).

We end this section by stating a fact that is proven in \cite{CFL}. Let us first recall that a \emph{scattered} linear ordering is a linear ordering that does not contain a copy of $\mathbb{Q}$. First let us recall that a classical representation theorem by Hausdorff provides a general expression of linearly ordered sets in terms of scattered sets (see \cite{rosenstein}, Theorem 4.9):

\begin{thm}[Hausdorff]
Any linear ordering \( L \) is a dense sum of scattered linear orderings; that is, there is a dense linear ordering \( L_* \) and a map \( h \) from \( L_* \) to scattered linear orderings such that $$ L = \sum \{ h(i) \mid i \in L_* \}. $$
\end{thm} 

We record here, without proof, a result from \cite{CFL}, obtained using the above classical theorem and describes the EOS for transitive homeomorphisms in more detail. Note that we choose to state the result for completeness, but nothing in the present paper depends on this fact.

\begin{fact}\label{long_t}
    Let $(X,f)$ be a compact dynamical system with $\text{card}(X)=\mathfrak{c}$. If $f$ is a transitive homeomorphism, then every countable  scattered ordering and the order-type of the rationals $\eta$ appear in $\Omega(X^2)$.
    
\vspace{0.3cm}

    More precisely:
    \begin{itemize}
        \item[1.]There exists a co-meagre set $\mathcal{S}\subseteq X^2$ such that, for every $(x,y)\in \mathcal{S}$, the family of orderings $\Omega(x,y)$ contains every scattered countable infinite ordering;
        \item[2.] $\Omega(x,y)$ contains the finite ordinal $k$ if and only if $x\,\mathcal{O}\,y$
        \item[3.]$\Omega(x,y)$ contains the order-type $\eta$ for every $x,y\in X$.
    \end{itemize}
\end{fact}
 The proof is quite long and proceeds by transfinite induction on the rank of the scattered order. We refer the reader to \cite{CFL} for the proof.
 

\bigskip

Now we want to analyze the order spectrum in terms of dual attractor/repeller pairs. For this, we need some preliminary results describing the behavior of nested chains on the gradient-like part of the system. We will use indeed the following result, which is stated (for flows) in \cite{Jacobs} (p. 2). To get the result for discrete iteration of maps it is enough to observe that $x\,\mathcal{C}\,y$ implies $f(x)\,\mathcal{C}\,y$ unless $f(x)=y$.
\begin{thm}\label{prop:Jacobs}
    Suppose that $x\,\mathcal{C}\,y$, Then exactly one of the 
    following alternatives holds:
    \begin{enumerate}
        \item[(1)] $y$ lies on the forward orbit of $x$, i.e.\ 
        $y=f^k(x)$ for some $k\in\mathbb{N}$;
        \item[(2)] for every $k\in\mathbb{N}$, we have 
        $f^k(x)\,\mathcal{C}\,y$.
    \end{enumerate}
\end{thm}

We also recall two well-known facts concerning attractors (see for instance \cite{kurka03}, p. 80-83.)

\begin{lem}\label{lem:omega-neighbourhood}
    Let $(X,f)$ be a compact dynamical system and $x\in X$. 
    If $V\subseteq X$ is an open set such that $\omega(x)\subseteq V$, then 
    there exists $N\in\mathbb{N}$ such that $f^n(x)\in V$ for every $n\geq N$.
\end{lem}

\begin{lem}\label{lem:chain-trap}
    Let $(X,f)$ be a compact dynamical system and let $A\subseteq X$ be an 
    attractor with inward set $U$, i.e.\ $U$ is closed, 
    $f(U)\subseteq \operatorname{int}(U)$, and 
    $A=\bigcap_{n\geq 0} f^n(U)$. Then there exists $\varepsilon_0>0$ such that
    every $\varepsilon_0$-chain $C:x_0,\dots,x_m$ with $x_0\in U$ satisfies 
    $x_i\in U$ for all $i=0,\dots,m$.
\end{lem}

We now prove that the chain relation, on the basin of an attractor (but outside the attractor), coincides in fact with the orbit relation.

\begin{thm}\label{prop:gradient-like-orbit}
    Let $(X,f)$ be a compact dynamical system and let $A$ be an attractor with 
    inward set $U$ and basin $B(A)$. Suppose $x,y\in B(A)\setminus A$ and 
    $x\,\mathcal{C}\,y$. Then there exists $k\in\mathbb{N}$ such that 
    $y=f^k(x)$.
\end{thm}

\begin{proof}
    By Theorem~\ref{prop:Jacobs}, since $x\,\mathcal{C}\,y$, either
    \begin{enumerate}
        \item[(1)] $y=f^k(x)$ for some $k\in\mathbb{N}$, or
        \item[(2)] for every $k\in\mathbb{N}$ we have $f^k(x)\,\mathcal{C}\,y$.
    \end{enumerate}
    If (1) holds we are done. Thus, for a contradiction, assume that (2) holds
    and $y$ is not on the forward orbit of $x$.

    Since $A$ is an attractor and 
    $y\notin A$, we can choose a closed inward set $U$ such that  $A\subseteq \operatorname{int}(U)$, $A=\cap_{k\to \infty} f^k(U)$ and $y\notin U$. 
    As $x\in B(A)$, we have $\omega(x)\subseteq A\subseteq \operatorname{int}(U)$. By 
    Lemma~\ref{lem:omega-neighbourhood}, there exists $N\in\mathbb{N}$ such that
    $$
        f^n(x)\in \operatorname{int}(U) \subseteq U \qquad\text{for every }n\geq N.
    $$

    Let $\varepsilon_0>0$ be given by Lemma~\ref{lem:chain-trap} for the inward 
    set $U$. Then any $\varepsilon_0$-chain starting in $U$ stays in $U$. In 
    particular, for each $n\geq N$ there can be no $\varepsilon_0$-chain from 
    $f^n(x)$ to $y$, because $f^n(x)\in U$ while $y\notin U$.

    This contradicts alternative~(2) of Theorem~\ref{prop:Jacobs}, which 
    asserts that $f^n(x)\,\mathcal{C}\,y$ for every $n\in\mathbb{N}$, i.e.\ for 
    every $\varepsilon>0$ (hence for $\varepsilon_0$ in particular) there should 
    exist an $\varepsilon$-chain from $f^n(x)$ to $y$. Therefore alternative~(2) 
    cannot occur, and we must be in case~(1). Hence there exists $k\in\mathbb{N}$ 
    with $y=f^k(x)$.
\end{proof}

\begin{thm}
    Let $x,y\in X$ be such that $x\,\mathcal{R}\,y$ and $(x,y)\notin\mathcal{O}$. Suppose that $\mathcal{O}(y)$ is a stable, attractive periodic orbit with period $K$. Then, $\Omega(x,y)=\{\omega,\,\omega+1,\dots,\,\omega+K-1\}$.
\end{thm}
\begin{proof}
    First of all, let us prove $\{\omega,\,\omega+1,\dots,\,\omega+K-1\}\subseteq\Omega(x,y)$. Let $j\in\{0,\dots,K-1\}$ and take $z\in\mathcal{O}(y)$ to be such that $f^j(z)=y$.
    Since $\mathcal{O}(y)$ is an attractor and $x\,\mathcal{R}\,y$, there exists a strictly increasing sequence of natural numbers $\{k_n\}_{n=1}^\infty$ such that $f^{k_n+1}(x)\xrightarrow{n\to\infty}z$ and $f^{k_n}(x)\neq z$. Then, in order to define a complete sequence of order-compatible nested chains $\{S_n\}_{n=1}^\infty$, it is enough to take $S_n$ as
    $$S_n:\,x,\,f(x),\dots,\,f^{k_n}(x),\,z,\,f(z),\dots,\,f^{j-1}(z)=y.$$
    
    Now, we will show that $\omega+j$, for $j\in\{0,\dots,\,K-1\}$, are the only possible order types in $\Omega(x,y)$.
    Let $\{S_n\}_{n=1}^\infty$ be a complete sequence of order-compatible nested chains from $x$ to $y$ for a sequence $\{\varepsilon_n\}_{n=1}^\infty$ and set $S=\bigcup_n\widehat{S_n}$. We want to show that the order type of the ordered limit set $(S',\le_\infty)$, where $S'=S\setminus\{x,y\}$, is $\omega+j$ for some $j\in\{0,\dots,\,K-1\}$. The set $S'$ decomposes in
    $S'=S^A\sqcup S^B,$
    where $S^A=S'\setminus\mathcal{O}(y)$ and $S^B=S'\cap\mathcal{O}(y)$.
    By Lemma \ref{lem:ONC-implies-chain}, we have that $x\,\mathcal{C}\,z$, for every $z\in S^A$, and, by Theorem \ref{prop:gradient-like-orbit}, we have $x\,\mathcal{O}\,z$. Thus $S^A\subseteq\mathcal{O}(x)$. On the other hand, it is not difficult to see that $\mathcal{O}(x)\subseteq S^A$.
    
    Indeed, suppose there exists $k\in\mathbb N$ such that $f^k(x)\in S^A$. We now show that $f^{k+1}(x)\in S^A$, too. Since $\omega(x)=\mathcal{O}(y)$, there exists an $\varepsilon>0$ such that
    $$B_\varepsilon(f^{k+1}(x))\cap(\mathcal{O}(x)\cup\mathcal{O}(y))=\{f^{k+1}(x)\}.$$
    Then, for every $n\in\mathbb N$ such that $f^k(x)\in\widehat{S_n}$ and $\varepsilon_n<\varepsilon$, the point $f^k(x)$ must be followed by $f^{k+1}(x)$ in the chain $S_n$ and, hence, $f^{k+1}(x)\in S^A$.
    Therefore, the order type of $(S^A,{\le_\infty}_|)$ is $\omega$.

    To have the statement, it is enough to note that $\operatorname{card}(S^B)<K$. In fact, otherwise the chains $S_n$ would contain cyclic sub-chains.
\end{proof}

We now use Theorem \ref{prop:gradient-like-orbit} to prove that, in the basin of an attractor, but outside the attractor itself, order-compatible nested chains can only be supported on full orbits exiting from one chain component and entering into the other one. Let us make precise the statement.

\begin{defi}
    We call a \emph{full orbit} of $f$ a bi-infinite sequence
    $(z_n)_{n=-\infty}^{\infty}$ of points in $X$ such that
    $f(z_n) = z_{n+1}$ for all $n\in\mathbb{Z}$. We call a \emph{negative semi-orbit} of $f$ a one-sided infinite sequence $(z_n)_{n=-\infty}^{-1}$ of points in $X$ such that $f(z_n) = z_{n+1}$ for all $n<-1$.
\end{defi}
\begin{lem}\label{lem:nested-contained-in-orbit}
    Let $A\subseteq X$ be an attractor and let $R:=X\setminus B(A)$ be its dual repeller. Let $x\in R$, $y\in A$ be such that there is a sequence of order-compatible nested chains $\{S_n\}_{n=1}^{\infty}$ from $x$ to $y$ for a vanishing sequence of positive real numbers $\{\varepsilon_n\}_{n=1}^\infty$. Then, setting
    $$S=\bigcup_n \widehat{S_n}\quad\text{and}\quad S'=S\setminus(A\cup R),$$
    one of the following holds:
\begin{enumerate}
    \item there exists a full orbit $\{z_n\}_{n=-\infty}^\infty$ such that $S'=\{z_n\}_{n=-\infty}^\infty$;
    \item there exists a negative semi-orbit $\{z_n\}_{n=-\infty}^{-1}$ such that $S'=\{z_n\}_{n=-\infty}^{-1}.$
\end{enumerate}  
\end{lem}

\begin{proof}
Since $R$ is invariant, we have $f(x)\in R$, so $f(x)\neq y$. Since $$\inf \{d(x',x''):x'\in R, x''\in A\}>0,$$ for $\varepsilon>0$ small enough there must be a point in $S$ belonging to $B(A)$, and thus $S'\neq\varnothing$. Then, pick some point $z_0\in S'$. If $S'=\{z_0\}$, it means that for every $\varepsilon>0$ there exists a point $w\in R$ such that $d(f(w),z_0)<\varepsilon$ which is in  contradiction with the fact that $R$ is a repeller (and thus closed in particular).

Let $w\in S'$ be arbitrary, $w\neq z_0$. Since $\{S_n\}_{n=1}^{\infty}$ is
a complete sequence of order-compatible nested chains, there exists $N$ such that for all $n\ge N$ both
$z_0$ and $w$ appear in $S_n$ and always in the same relative order.
Thus, exactly one of the following holds:

\begin{itemize}
    \item[(a)] $z_0$ appears before $w$ eventually;
    \item[(b)] $w$ appears before $z_0$ eventually.
\end{itemize}

In case (a), Lemma~\ref{lem:ONC-implies-chain} gives
$z_0\,\mathcal{C}\,w$. Since $z_0,w\in B(A)\setminus A$, by Theorem~\ref{prop:gradient-like-orbit}, we have
$$w = f^{k}(z_0)$$
for some $k>0$.

In case (b), Lemma~\ref{lem:ONC-implies-chain} gives $w\,\mathcal{C}\,z_0$.
Again, arguing as above and using Theorem~\ref{prop:gradient-like-orbit}, we have
$$z_0 = f^{k}(w)$$
for some $k> 0$.

Summarizing, for every $w\in S'$ there exists an integer
$k(w)> 0$ such that either
$$
   w = f^{k(w)}(z_0)
   \quad\text{or}\quad
   z_0 = f^{k(w)}(w).
$$
Hence we have: 

\begin{equation}\label{chain_orbit_}
\forall\  z,w\in S', \ z\,\mathcal{O}\,w \iff z\le_\infty w.
\end{equation}

We now prove that every $z\in S'$ is the exact $f$-image of some point in $S'$, i.e.\ there exists $w\in S'$ with $f(w)=z$. In fact, if there was a point $z\in S'$ that is not the exact $f$-image of any point $w\in S'$, by \eqref{chain_orbit_} we would have that $z$ is the minimal element in $(S',{\le_\infty}_{|S'})$, which implies that $S'\subseteq\mathcal{O}(z)\cup\{z\}$ and in general $S\subseteq R\cup\mathcal{O}(z)\cup\{z\}\cup A$.

Since $z\in B(A)\setminus A$, $R$ is closed and since $\omega(z)\subseteq A$, there exists an $\varepsilon>0$ such that
$$B_\varepsilon(z)\cap(R\cup\mathcal{O}(z)\cup A)=\varnothing.$$
Now, for every $n\in\mathbb N$ such that $\varepsilon_n<\varepsilon$, we have $z\notin\widehat{S_n}$, because, for every $n\in\mathbb N$, we have $\widehat{S_n}\setminus\{z\}\subseteq R\cup\mathcal{O}(z)\cup A$ and $f(R\cup\mathcal{O}(z)\cup A)\subseteq R\cup\mathcal{O}(z)\cup A$. Thus, we found a contradiction and, therefore, every element $w$ in $S'$ has a predecessor with respect to $\le_\infty$ in $f^{-1}(w)\cap S'$.

We first show that this predecessor is unique. Suppose there exist $p,q\in S'$ with $p\neq q$ and
$$
f(p)=f(q)=z.
$$
Then
$$
\mathcal{O}(p)=\{f(p)=z,f(z),f^2(z),\dots\}\qquad ,\qquad 
\mathcal{O}(q)=\{f(q)=z,f(z),f^2(z),\dots\}.
$$
In particular,
$$
q\notin\mathcal{O}(p)\qquad ,\qquad p\notin\mathcal{O}(q),
$$
so $p$ and $q$ are not orbit-related, because otherwise $z$ would be a periodic point in $B(A)\setminus A$. Thus, we found a contradiction, because we proved that, for any two points $z,w\in S'$, we have $z\,\mathcal{O}\,w$ or $w\,\mathcal{O}\,z$. Hence each $z\in S'$ has exactly one predecessor in $S'$. It follows that $f|_{S'}$ is injective on $S'$.

Therefore, we have
$$
z\le_\infty w \quad\Longleftrightarrow\quad z=w \text{ or } w\in\mathcal{O}(z).
$$
For any $z,w\in S'$ we have $z\le_\infty w$ or $w\le_\infty z$, so the relation is total. It is obviously also transitive and, since there are no periodic points in $S'\subseteq B(A)\setminus A$, it is antisymmetric, too.  Thus $\le_\infty$ is a linear order on $S'$. Denote
$$
z<_\infty w \quad\Longleftrightarrow\quad z\le_\infty w,\;z\neq w.
$$

Let $w,z\in S'$ with $f(w)=z$. We claim that there is no $t\in S'$ with
$$
w<_\infty t<_\infty z.
$$
Assume by contradiction that such $t$ exists. By linearity of the order, $w<_\infty t$ implies $t\in\mathcal{O}(w)$, and $t<_\infty z$ implies $z\in\mathcal{O}(t)$. Thus there exist $m,k>0$ such that
$$
t=f^m(w),\qquad z=f^k(t).
$$
On the other hand we also know $z=f(w)$, hence
$$
f(w)=f^k(t)=f^k(f^m(w))=f^{m+k}(w).
$$
Since $f_{|S'}$ is injective, we deduce $f^{m+k-1}(w)=w$, and so $w$ is a periodic point living in $B(A)\setminus A$, which is absurd.

Thus, no such $t$ exists, and $z$ is the immediate successor of $w$ with respect  to the order $\le_\infty$. In particular, for each $z\in S'$, there is a unique $w\in S'$ such that $f(w)=z$, and this $w$ is the immediate predecessor of $z$.

We have shown that $(S',\le_\infty)$ is a linearly ordered set such that:
\begin{itemize}
\item every element has a unique immediate predecessor (coming from the unique $w$ such that $f(w)=z$),
\item there is no minimal element.
\end{itemize}
There are two cases.

\smallskip

\textit{Case 1: there is no maximal element in $(S',\le_\infty)$.} Then, the ordered set $(S',\le_\infty)$ is such that every element $z\in S'$ has a unique immediate successor $w\in S'$; moreover, we have $w=f(z)$.
Then, the ordered set $(S',\le_\infty)$ is order-isomorphic to $\mathbb{Z}$. Hence, there exists a bijection
$$
\varphi:\mathbb{Z}\to S',\qquad k\mapsto z_k,
$$
such that
$$
k<\ell \;\Longleftrightarrow\; z_k<_\infty z_\ell,
$$
and the immediate successor of $z_k$ in $S'$ is $z_{k+1}$. By construction of the order, this successor relation coincides with the dynamics, so
$$
f(z_k)=z_{k+1}\quad\text{for all }k\in\mathbb{Z},
$$
and
$$
S'=\{z_k:k\in\mathbb{Z}\}.
$$

\smallskip

\textit{Case 2: there is a maximal element.} Then, there exists a bijection $\psi$ between the negative integers $\mathbb Z\setminus\mathbb N_0$ and $S'$
$$\psi:k\mapsto z_k$$
such that
$$k<\ell\iff z_k<_\infty z_\ell.$$
Also in this case, by the construction of the order, we have
$$f(z_k)=z_{k+1}\quad\text{for all }k\in\mathbb Z\setminus\mathbb N_0,$$
and so
$$S'=\{z_k\mid k\in\mathbb Z\setminus\mathbb N_0\}.$$

\end{proof}

\begin{thm}
Let $A$ be an attractor and $R$ its dual repeller. Let $x\in R$ and $y\in A$. Then every $\tau\in \Omega(x,y)$ admits the decomposition: $$\tau = \beta + \eta + \beta'$$
where: 
\begin{itemize} 
\item $\beta, \beta'$ correspond to the order type of nested chains supported on $R$ and $A$ respectively;
\item $\eta$ is equal to $\omega^*$ or $\zeta$ and corresponds to the order-type of nested chains supported on $B(A)\setminus A$.
\end{itemize}

\end{thm}

\begin{proof}
If $\Omega(x,y)= \emptyset$ there is nothing to prove. Assume thus $\Omega(x,y)\neq \emptyset$.

Let $\{S_n\}_{n=1}^\infty$ be a complete sequence of order-compatible nested chains from $x$ to $y$ such that the ordered limit set $(S\setminus\{x,y\},\le_\infty)$, where $S=\bigcup_n\widehat{S_n}$, has order-type $\tau$. Then, the set $S$ decomposes in
$$S=S_R\sqcup S'\sqcup S_A,$$
where
$$S_R=S\cap R\quad,\quad S_A=S\cap A\quad\text{and}\quad S'=S\setminus(A\cup R).$$
By Lemma \ref{lem:nested-contained-in-orbit}, there exists either a full orbit $O=\{z_k\}_{k=-\infty}^\infty$ or a negative semi-orbit $O=\{z_k\}_{k=-\infty}^{-1}$ such that $S'=O$. Then, clearly, the order type of $(S',\le^* _\infty)$, where $\le^* _\infty$ is the order relation $\le_\infty$ restricted to $S'$, is, respectively, $\zeta$ or $\omega^*$.

Thus, to prove the thesis, it is enough to define $\beta$ as the order type of $\left(S_R\setminus\{x\},\le^R_\infty\right)$, where $\le^R_\infty$ is the order relation $\le_\infty$ restricted to $S_R\setminus\{x\}$, and to define $\beta'$ as the order type of $\left(S_A\setminus\{y\}),\le_\infty^A\right)$, where $\le_\infty^A$ is the order relation $\le_\infty$ restricted to $S_A\setminus\{y\}$.
\end{proof}

\section{Refinement of Conley decomposition and prolongational hierarchy}
\noindent For a compact system $(X,f)$, the EOS $\Omega$ witnesses the chain recurrent part of the system, in the sense that, if $\Omega(x,y)\neq\varnothing\neq\Omega(y,x)$, then $x,y$ lie in the same chain component. Let us see how this basic fact can be refined.

First of all, if $K$ is a chain component, the spectra observed inside $K$ are intrinsic, as we prove in the following result, to the subsystem $(K,f|_K)$, independent of the ambient system, and provide a canonical decoration of each Conley component. 

\begin{thm}
    Let $(X,f)$ be a compact dynamical systems, let $K\subseteq X$ be a chain component and consider the dynamical system $(K,f|_{K})$. Then,
    $$\big(\Omega\big)_{f_{|_K}}=\big(\Omega_{f}\big)_{|_{K^2}}$$
\end{thm}
\begin{proof}
    Trivially, for every $(x,y)\in K^2$, we have $\Omega_{f_{|K}}(x,y)\subseteq\Omega_f(x,y)$, because any $\varepsilon$-chain from $x$ to $y$ within $K$ is an $\varepsilon$-chain in the dynamical system $(X,f)$, too.
    
    To prove the other inclusion,
    we show that, given $(x,y)\in K^2$, any order-compatible nested sequence $\{S_n\}_{n=1}^\infty$ from $x$ to $y$ is such that $\widehat{S_n}\subseteq K$ for every $n\in\mathbb N$. From this fact, it follows that for every ordering $\beta\in\Omega_f(x,y)$, we have $\beta\in\Omega_{f_{|K}}(x,y)$, too. Moreover, by Lemma \ref{lem:ONC-implies-chain}, for every $z\in\bigcup_n\widehat{S_n}$, we have $x\,\mathcal{C}\,z$ and $z\,\mathcal{C}\,y$, and since $x$ and $y$ belong to the same chain component $K$, we have $z\in K$, too. 
\end{proof}

In the next two results, we will see that the limit order $(S',\leq_\infty)$ attached to a sequence of order-compatible
nested chains from $x$ to $y$ (where $S'=\left(\cup_n\widehat{S_n}\right)\setminus\{x,y\}$) induces a  refinement of Conley’s partial order.
Indeed, the components that are
actually visited by the chains appear as contiguous ``blocks'' in $(S',\le_\infty)$, and
the quotient that collapses each such block
to a single point recovers the Conley order on the corresponding set of
components.
The next two results make this precise: first we show that the intersection of
each chain component with $S'$ is convex in $(S',\le_\infty)$, and then that the
induced Conley order on the corresponding set of components is in fact linear.

 \begin{thm}\label{convex_}
    Let $(X,f)$ be a compact dynamical system and let $x,y\in X$ be such that $x\,\mathcal{C}\,y$. Let $\{S_n\}_{n=1}^\infty$ be a complete sequence of order-compatible nested chains and set $S=\bigcup_n\widehat{S_n}$ and $S'=S\setminus\{x,y\}$. Then, for every chain component $K\subseteq X$, the set $K\cap S'$ is convex in the ordered limit set $(S',\le_\infty)$.
\end{thm}
\begin{proof}
    Pick $x_1,x_2\in K\cap S'$ such that $x_1\le_\infty x_2$ and take $x_3\in S'$ such that $x_1\le_\infty x_3\le_\infty x_2$. Then, we need to prove that $x_3\in K$. In fact, by Lemma \ref{lem:ONC-implies-chain}, we have $x_1\,\mathcal{C}\,x_3$  and $x_3\,\mathcal{C}\,x_2$ and, since $x_1,x_2\in K$, we have also $x_2\,\mathcal{C}\,x_1$. Therefore, $x_3\,\mathcal{C}\,x_1$, because $\mathcal{C}$ is a transitive relation, and thus, $x_3\in K$.
\end{proof}

\begin{thm}\label{conley_}
    Let $(X,f)$ be a compact dynamical system and let $x,y\in X$ be such that $x\,\mathcal{C}\,y$. Let $\{S_n\}_{n=1}^\infty$ be a complete sequence of order-compatible nested chains and set $S=\bigcup_n\widehat{S_n}$ and $S'=S\setminus\{x,y\}$. Let $H$ be the set of all the chain components $K\subseteq X$ such that $K\cap S\neq\varnothing$.
    Then, $(H,\le_\mathrm{Conley})$ is a linearly ordered set.
\end{thm}

\begin{proof}
    It is known that  $(H,\le_{\mathrm{Conley}})$ is a partially ordered set. It remains to prove that for every $K, K'\in H$ one of the following holds: $K\le_{\mathrm{Conley}}K'$ or $K'\le_{\mathrm{Conley}}K$.
    We distinguish two cases.
    If $x=y$, it means that there is a chain component $K$ such that $S\subseteq K$, so $H=\{K\}$ and so there is nothing to prove. In fact, let $z\in S\setminus\{x\}$ and let $S_n:\,x_0^{(n)}=x,\,x_1^{(n)},\dots,\,x_{m_n}^{(n)}=x$. Thus, there exists a sequence of points $\{x_{i_n}^{(n)}\}_{n=1}^\infty$, with $i_n\in\{1,\dots,m_n-1\}$, such that $x_{i_n}^{(n)}=z$. Then
    $$\{S_n^{x\to z}:\,x_0^{(n)},\,x_1^{(n)},\dots,\,x_{i_n}^{(n)}\}_{n=1}^\infty\quad\text{and}\quad\{S_n^{z\to x}:\,x_{i_n}^{(n)},\,x_{i_n+1}^{(n)},\dots,\,x_{m_n}^{(n)}\}_{n=1}^\infty$$
    are two complete sequences of chains respectively from $x$ to $z$ and from $z$ to $x$, and therefore $z\in K$.
    
    If $x\ne y$, we can have that $S\subseteq K$ and so $H=\{K\}$ and the thesis is trivially valid. Instead, if $x\ne y$ and there are at least two different chain components $K$ and $K'$ such that both $K\cap S$ and $K'\cap S$ are not empty, we
    take two distinct points $z\in K\cap S$ and $w\in K'\cap S$. Suppose that, for all sufficiently large $n$, $z$ appears before $w$ in $S_n$. Then, by Lemma \ref{lem:ONC-implies-chain}, we have that $z\,\mathcal{C}\,w$ and thus $K\le_{\mathrm{Conley}}K'$. Similarly, if $w$ appears before $z$, we can conclude that $K'\le_{\mathrm{Conley}} K$. 
\end{proof}
Note that in the last two results we did not use order-compatibility, but just weak order-compatibility for the nested, acyclic chains, which is always achievable up to subsequences.

\bigskip

The most important transfinite structure used to describe recurrence is probably the prolongational hierarchy $J_\alpha(x)$ introduced by Auslander (see \cite{Auslander}). We recall that the prolongational set of $x$ of order $\alpha$ is defined, in the discrete-time case, as (see \cite{Viscovini}):
\begin{equation}
J_1(x)=\{y\in X\,|\,x\,\mathcal{N}\,y\}\quad,\quad
   J_\alpha(x) = \bigcap_{\varepsilon > 0} \overline{ \bigcup_{n=1}^\infty \bigcup_{\beta < \alpha} (J_\beta)^n(B_\varepsilon(x)) },
\end{equation}
where $(J_\beta)^1=J_\beta$ and $(J_\beta)^n=(J_\beta)((J_\beta)^{n-1})$.

We want to show some simple examples in which the EOS provides a finer description of recurrence properties than the prolongational sets. The fact that any emergent order in $\Omega(x,y)$ is ``anchored" to a precise sequence of ordinately nested chains makes it possible to describe the recurrence between a certain ordered pair of points in distinct ways, each linked to a different sequence of chains, whereas different transfinite levels of prolongational sets are strictly nested: if $y\in J_\alpha(x)$ then it automatically belongs to $J_\beta(x)$ for every countable $\beta>\alpha$. For instance, if a point $x$ is in recurrence relation with $y\notin\mathcal{O}(x)$, then clearly $y\in J_1(x)$ (and thus to $J_\alpha$ for every $\alpha$). This may happen, to mention two extreme cases, either if: 
\begin{itemize}
    \item  $x$ is in the basin of the attracting fixed point $y$, 
    
    or if
    \item $x$ is a transitive point.
\end{itemize}
From the point of view of the EOS, we have in the two cases:
\begin{itemize}
    \item $\Omega(x,y)$ coincides with the unique order-type $\{\omega\}$;
    \item $\Omega(x,y)$ contains typically many other order-types; for instance, $\omega+k$ is in $\Omega(x,y)$ for every $k\in\mathbb{N}$.
\end{itemize}

More generally, prolongational sets only ``count" the needed $\varepsilon$-\textit{corrections}, while emergent orders directly count \textit{iterations}, and this can provide a finer description of recurrence. We now discuss some examples of dynamical systems $(X,f)$ in which we have three points $x_1,x_2,x_3\in X$ where we cannot distinguish $x_2$ and $x_3$ from the point of view of prolongational set of $x_1$, while the EOS discriminates the recurrence properties between $x_1$ and $x_2$ and between $x_1$ and $x_3$.

\begin{thm}
    There exists a dynamical system $(X,f)$ and three points $x_1,x_2,x_3\in X$ such that
    $$x_3\in J_\alpha(x_1)\quad,\quad x_3\in J_\alpha(x_2)$$
    for some ordinal $\alpha=\min\{\beta\mid x_3\in J_\beta(x_1)\}=\min\{\beta\mid x_3\in J_\beta(x_2)\}$, and
    $$\Omega(x_1,x_3)\neq\Omega(x_2,x_3).$$
\end{thm}
\begin{proof}
    We will prove the statement by exhibiting an instance of the phenomenon.

    Consider the dynamical system
    $(X,f)$ (a commonly used example of a countable collection of chain components in a continuum, see Fig.\ref{fig_1}), where: 
    \begin{equation*}
        X=[0,1]\quad,\quad f(x)=2^{-\lfloor \log_2 x \rfloor}\bigl(x-2^{\lfloor \log_2 x \rfloor}\bigr)^2+2^{\lfloor \log_2 x \rfloor},
        \label{interval_map}
    \end{equation*}
    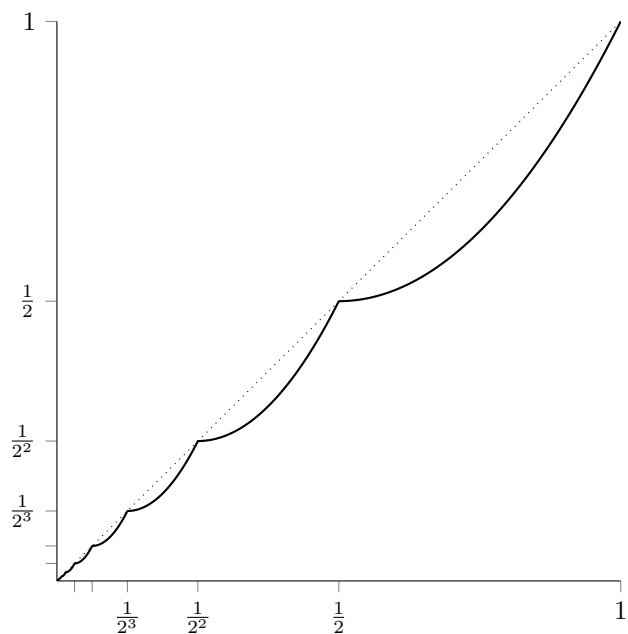
\begin{figure}[ht!] 
        \centering
        \begin{tikzpicture}
            \begin{axis}[
                width=9cm, height=9cm,
                axis equal,
                xmin=0, xmax=1,
                ymin=0, ymax=1,
                domain=0.0001:1,
                samples=700,
                axis lines=middle,
                axis line style={-},
                xtick={1,0.5,0.25,0.125,0.0625,0.03125},
                xticklabels={$1$,$\frac{1}{2}$,$\frac{1}{2^2}$,$\frac{1}{2^3}$},
                ytick={1,0.5,0.25,0.125,0.0625,0.03125},
                yticklabels={$1$,$\frac{1}{2}$,$\frac{1}{2^2}$,$\frac{1}{2^3}$},
                tick align=outside,
                ]
                \addplot[thick, black] 
                { 2^(-floor(ln(x)/ln(2))) * ( x - 2^(floor(ln(x)/ln(2))) )^2 + 2^(floor(ln(x)/ln(2))) };
                \addplot[dotted, black] 
                { x };
            \end{axis}
        \end{tikzpicture}
        \caption{Plot of the interval map defined in Eq.\eqref{interval_map}.}
        \label{fig_1}
    \end{figure}
    Then, setting $x_1=1$ and $x_3=0$ and picking $x_2\in(1/2,1)$, we have
    $$x_3\in J_2(x_1)\setminus J_1(x_1)\quad\text{and}\quad x_3\in J_2(x_2)\setminus J_1(x_2),$$
    and
    $$\Omega(x_1,x_3)\neq\Omega(x_2,x_3),$$
    because $x_3\in[\zeta\cdot\omega](x_1)\setminus[\omega+\zeta\cdot\omega](x_1)$ and $x_3\in[\omega+\zeta\cdot\omega](x_2)\setminus[\zeta\cdot\omega](x_2)$.  
\end{proof}

 Let us now show another instance of the refinement we are describing, this time using the classical Denjoy homeomorphism (\cite{denjoy}).

 \begin{defi}[Denjoy circle homeomorphism]
 $$
 I_n:=f^n(I)\qquad(n\in\mathbb Z).
 $$
 \end{defi}

 \begin{fact}
 Let $f$ be a Denjoy homeomorphism and $I$ a wandering interval, with $I_n=f^n(I)$. Then:
 \begin{enumerate}
 \item There exists a continuous, surjective, degree–one monotone map $h:\mathbb S^1\to\mathbb S^1$ with
 $$ h\circ f \;=\; R_{\rho(f)}\circ h, $$
 where $R_{\rho(f)}$ is the rigid rotation by angle $\rho(f)$. The map $h$ collapses each component of
 $\bigcup_{n\in\mathbb Z} I_n$ to a point; hence $f$ is semi–conjugate but not conjugate to $R_{\rho(f)}$.
 \item The closed set
 $$ K \;:=\; \mathbb S^1 \setminus \bigcup_{n\in\mathbb Z} I_n $$
 is a perfect, totally disconnected, nowhere–dense (Cantor) $f$–invariant set, and $f|_K$ is minimal. All wandering intervals are exactly the family $\{I_n\}_{n=-\infty}^{\infty}$, and their lengths tend to $0$.
 \end{enumerate}
 \end{fact}

 \begin{thm}\label{denjo_}
     There exists a dynamical system $(X,f)$ and three points $x_1,x_2,x_3\in X$ such that
    $$x_2,x_3\in J_\alpha(x_1),$$
     for some ordinal $\alpha=\min\{\beta\mid x_2\in J_\beta(x_1)\}=\{\beta\mid x_3\in J_\beta(x_1)\}$, and
     $$\Omega(x_1,x_2)\neq\Omega(x_1,x_3).$$

 \end{thm}
 \begin{proof}
     We will prove the statement by exhibiting two instances of the phenomenon. 
     \begin{enumerate}
        \item Consider the
        dynamical system $(X,f)$, where $X=[-1,1]$ and
        $$
        f(x)=
        \begin{cases}
            2^{-\lfloor \log_2 x \rfloor}\bigl(x-2^{\lfloor \log_2 x \rfloor}\bigr)^2+2^{\lfloor \log_2 x \rfloor} &\text{if }0\le x\le1\\
            (x+1)^2-1&\text{if }-1\le x\le0
        \end{cases}.
        $$
        Then, setting $x_1=1$ and $x_2=-1$ and picking $x_3\in(-1,0)$, we have
        $$x_2,x_3\in J_3(x_1)\setminus J_2(x_1),$$
        and
        $$\zeta\cdot\omega+\zeta\in\Omega(x_1,x_2)\setminus\Omega(x_1,x_3)\quad\text{and}\quad\zeta\cdot\omega+\omega^*\in\Omega(x_1,x_3)\setminus\Omega(x_1,x_2),$$
        and thus, $\Omega(x_1,x_2)\neq\Omega(x_1,x_3)$.

        \item Let $f:\mathbb S^1\to\mathbb S^1$ be a Denjoy homeomorphism and consider the dynamical system $(\mathbb S^1,f)$. Pick $x_1,x_2\in K$ and $x_3\in I_0$. 
        Then, since $\{I_{-h}\}_{h\in \mathbb{N}}$ is dense, in every neighbourhood of $x_1$ we can find a point $z$ such that $z\,\mathcal{O}\,x_3$, thus $x_1\,\mathcal{N}\,x_3$ and $x_3\in J_1(x_1)$; of course we have also $x_2\in J_1(x_1)$, so $x_1$ and $x_2$ are indistinguishable based on their prolongational relation with $x_3$. Obeserve now that no point in the circle has $x_3$ in its $\omega$-limit set (no forward orbit accumulates at $x_3$), because $x_3$ belongs to the wandering set, so $\zeta\notin \Omega(x_1,x_3)$. In fact, if $\zeta$ was in $\Omega(x,y)$, then, given a sequence of order-compatible nested chains $\{S_n\}_{n=1}^\infty$ and named $S=\cup_n\widehat{S_n}$, the set $S\setminus\{x,y\}$ would consist of a full-orbit $(z_n)_{n=-\infty}^\infty$, because every element in $\zeta$ has both predecessor and successor. But, this would imply that $\mathcal{O}(z_0)$ accumulates at $x_3$ which is impossible.
        On the other hand, it is not difficult to see that $\Omega(x_1,x_2)$ contains $\zeta$.
    \end{enumerate}
    
\end{proof}

\begin{rem}
The previous example shows that the implication in Theorem \ref{zetaN} cannot be inverted in general.
\end{rem}
\vspace{0.3cm}
Note that $\Omega$ elementarily discriminates between the two dynamical systems $(\mathbb S^1,f)$, where $f$ is the Denjoy homeomorphism, and $(\mathbb S^1,R_\alpha)$, where $R_\alpha$ is the irrational rotation. For instance, $\zeta\in\Omega_{R_\alpha} (x,x)$ for every $x\in\mathbb S^1$, while $\zeta\notin\Omega_f(x,x)$ if $x$ is in a wandering interval $I$.

\bigskip

 \section{Declarations}

 The authors declare that they have no conflicts of interest.

 This research did not receive any specific grant from funding agencies in the public, commercial, or not-for-profit sectors.

\end{document}